\documentclass[11pt]{article}
\usepackage[T1]{fontenc}
\usepackage{euscript}
\usepackage{amsmath}
\usepackage{amssymb}
\usepackage{latexsym}
\usepackage{fancyhdr}

\usepackage{diagrams}
\include{amsfonts}

\newarrow{equal}=====
\newarrow{to}---->
\newarrow{onto}----{>>}
\newarrow{embed}>--->
\newarrow{TeXonto}{-}{-}{-}{-}{->>}
\newarrow{TeXto}----{->}
\newarrow{TeXembed}{>}---{->}
\newarrow{TeXembedd}{>}---{->>}

\newcounter{llistadepth}

\newenvironment{manlist}[1]{\addtocounter{llistadepth}{1}
      \edef\llistacontador{llista\romannumeral\the\value{llistadepth}}
      \list{{\rm ({#1{\llistacontador}})}}{\usecounter{\llistacontador}
      \def\makelabel##1{\hss\llap{##1}}
      \itemsep=2pt\parsep=0pt\topsep=3pt plus 1pt minus 1 pt}}{\endlist
      \addtocounter{llistadepth}{-1}}

\renewenvironment{enumerate}{\begin{manlist}{\roman}}{\end{manlist}}

\newtheorem{df}{Definition}[section]      
\newtheorem{thm}[df]{Theorem}             
\newtheorem{prop}[df]{Proposition}
\newtheorem{cor}[df]{Corollary}

\newtheorem{ex}[df]{Example}

\newcommand{\pf}{\noindent{\sc Proof.}\ }

\newcommand{\boom}{\quad\lower3pt\hbox{\vrule height1.1ex width .9ex depth -.2ex}
                    \vskip9pt}

\renewcommand{\mathcal}[1]{\EuScript{#1}}

\let\phi=\varphi
\let\Ga=\Gamma

\let\Ups=\Upsilon

\newcommand{\ad}{\mathop{\rm ad}\nolimits}
\newcommand{\Ad}{\mathop{\rm Ad}\nolimits}
\def\Alt{{\rm Alt}}

\newcommand{\CDO}{\mathop{\rm CDO}}
\def\Der{\mathop{\rm Der}}

\def\dim{\mathop{\rm dim}}

\def\Exp{\mathop{\rm Exp}\,}

\def\Hom{\mathop{\rm Hom}}
\def\id{{\rm id}} 

\def\Int{\mathop{\rm Int}}

\def\rk{\mathop{\rm rk}}

\def\chigh{{\raise1.5pt\hbox{$\chi$}}}

\newcommand{\gog}{\mathfrak{g}}
\newcommand{\hoh}{\mathfrak{h}}

\newfont{\pointwise}{lcircle10 scaled 500}
\newcommand{\pwise}{{\mbox {\pointwise \symbol{118}}}}
\def\ptwise{\hskip.02in\raise2pt\hbox{$\pwise$}}

\newfont{\numbers}{bbm9 scaled 1200}
\newcommand{\reals}{{\mbox {\numbers R}}}
\newcommand{\integers}{{\mbox {\numbers Z}}}

\newcommand{\co}{\colon\thinspace}                

\newcommand{\st}{\ \vert\ }

\let\isom=\cong

\def\act{\mathbin{\hbox{$<\kern-.4em\mapstochar\kern.4em$}}}
\def\ract{\mathbin{\hbox{$\mapstochar\kern-.3em>$}}}

\let\Bar=\overline
\let\Hat=\widehat
\let\Tilde=\widetilde

\renewcommand{\to}{\longrightarrow}

\def\gpd{\,\lower1pt\hbox{$\longrightarrow$}\hskip-.24in\raise2pt
             \hbox{$\longrightarrow$}\,}

\newcommand{\surj}{-\!\!\!-\!\!\!-\!\!\!\gg}
\newcommand{\inj}{>\!\!\!-\!\!\!-\!\!\!-\!\!\!>}

\begin{document}
\title{{\bf On the connection theory of extensions of transitive Lie groupoids}
\thanks{Research supported by the Greek State Scholarships Foundation}}

\author{Iakovos Androulidakis\\
        C.N.R.S. - Universit\'{e} Pierre et Marie Curie (Paris VI)\\
        Institut de Math\'{e}matiques de Jussieu - UMR 7586\\
        175, Rue du Chevaleret\\
        F-75013 PARIS\\
        France\\
        {\sf androulidakis@math.jussieu.fr}}

\date{{\sf \today}}

\maketitle

\abstract{Due to a result by Mackenzie, extensions of transitive Lie groupoids are equivalent to certain Lie groupoids which admit an action of a Lie group. This paper is a treatment of the equivariant connection theory and holonomy of such groupoids, and shows that such connections give rise to the transition data necessary for the classification of their respective Lie algebroids.}


\section{Introduction}

A groupoid, in the categorical sense, is a category where every arrow is invertible. Lie groupoids are those ones whose object and arrow spaces are manifolds, and the structure maps (source $s$, target $t$, multiplication and inversion) are smooth, plus $t$ and $s$ are submersions. A Lie group can be considered as a Lie groupoid (with a single point as the object space), but in general Lie groupoids are inherently noncommutative, hence their extensive use in differential geometry. In the same fashion, a bundle of Lie groups is also a Lie groupoid, with its projection playing the role of both the source and target maps. In this paper we focus on the transitive case, namely those Lie groupoids $\Omega \gpd M$ such that the anchor map $(t,s) : \Omega \rightarrow M \times M$ is a surjective submersion. 

Given a transitive Lie groupoid $\Omega \gpd M$, an extension of this is a pair of Lie groupoid morphisms
\begin{eqnarray}\label{gpdextn2}
F \stackrel{\iota}{\inj} \Omega \stackrel{\phi}{\surj} \Phi
\end{eqnarray}
where $\Phi \gpd M$ is a transitive Lie groupoid and $F \rightarrow M$ is a bundle of Lie groups. Much of the work of Mackenzie concerns such extensions, and for a general Lie groupoid there always exists a foliation whose restriction to every leaf is such an extension. When dealing with such an extension, it is quite difficult to keep track of the transition functions of all the groupoids involved. Mackenzie in \cite{Mackenzie:Cahiers}, managed to reformulate such an extension to a single (transitive) Lie groupoid over the total space of a principal bundle, plus a group action. These are called PBG-groupoids, and an account of this correspondence is given in Section 1. It is therefore reasonable to expect that extensions of transitive Lie groupoids are classified by the transition functions of their corresponding PBG-groupoids, but it is necessary to keep track of the group action as well, and it is not certain that there exist transition functions which enjoy some kind of equivariance.

This is where connections are necessary. The transversals of an extension of transitive Lie groupoids correspond to connections in the respective PBG-groupoid which are suitably equivariant, and these are the connections whose infinitesimal theory and holonomy are given in this paper. The term {\it isometablic} is used for these connections, instead of equivariant, in order to highlight the non-standard nature of equivariance for the transition functions that they induce. 

The classification of PBG-groupoids (extensions of transitive Lie groupoids) is a problem of a different order presented in \cite{class}. In this paper we focus on the infinitesimal level, namely the classification of transitive PBG-algebroids. It is shown in \cite{LGLADG} that every transitive Lie algebroid $A$ (i.e. an extension of $TM$ by a Lie algebra bundle $L$) locally can be written as $TU_{i} \oplus L_{U_{i}}$, and in \cite{LGLADG} it is shown that such Lie algebroids are classified by pairs $(\chi,\alpha) = \{(\chi_{ij},\alpha_{ij})\}_{i,j \in I}$, where $\alpha_{ij} : U_{ij} \rightarrow Aut(\hoh)$ are the transition functions of the sections of $L$ and $\chi_{ij} : TU_{ij} \rightarrow U_{ij} \times \hoh$ are certain differential 1-forms, and the two of them satisfy suitable compatibility conditions. Here $\hoh$ is the fiber type of $L$ and $\{U_{i}\}_{i \in I}$ is a simple open cover of $M$. This data arises from the fact that locally there exist flat connections, namely Lie algebroid morphisms $TU_{i} \rightarrow A_{U_{i}}$. In this paper we show that locally PBG-algebroids have flat isometablic connections, and this gives rise to transition data which is suitably equivariant, modulo a simple connectivity assumption. Namely, we prove that:
\begin{thm}\label{0.1}
Let $P(M,G)$ be a principal bundle whose structure group $G$ is simply connected and $\{P_{i} \equiv U_{i} \times G\}_{i \in I}$ a local trivialisation of this bundle. If $A$ is a PBG-algebroid over this bundle then there exist transition data $(\chi,\alpha)$ such that 
\begin{enumerate}
\item $\chi_{ij}(Xg) = \chi_{ij}(X)g$ and
\item $\alpha_{ij}(ug) = \alpha_{ij}(u)g$
\end{enumerate}
for all $X \in TP_{ij}$, $u \in P_{ij}$ and $g \in G$.
\end{thm}

When $A$ is integrable though, the simple connectivity assumption is no longer necessary. It will be shown that:
\begin{thm}\label{0.2}
For any integrable PBG-algebroid $A$ over a principal bundle $P(M,G)$ there exist transition data $(\chi,\alpha)$ which satisfy $(i)$ and $(ii)$ of \ref{0.1}.
\end{thm}
Note that the in order to prove the above result, the equivariance of the sections of the PBG-groupoid which integrates $A$ is investigated, and this is an important step to the classification of PBG-groupoids. 

Finally, we give an account of the holonomy of isometablic connections which leads to the proof of an Ambrose-Singer theorem for isometablic connections. The treatment of holonomy presented here follows the fashion of \cite{LGLADG}, i.e. we use the equivariant version of deformable sections. This method provides more information for the relation of isometablic connections of PBG-Lie algebra bundles with their sections.

The structure of this paper is as follows: Sections 2 and 3 are an account of the correspondence between extensions and PBG structures, in both the groupoid and the algebroid level. A number of examples is included, most of which are used in the proofs later in the paper. Section 4 gives the connection theory of PBG-groupoids and PBG-Lie algebra bundles. In Section 5 we give the proofs of \ref{0.1} and \ref{0.2}. Section 6 is an account of the holonomy of isometablic connections in the fashion of \cite{LGLADG}, and Section 7 gives the proof of the Ambrose-Singer theorem for isometablic connections. 

\section*{Acknowledgments}
I would like to thank Kirill Mackenzie for the encouragement to write an 
account of connection theory in the PBG setting, as well as his useful comments. Also many thanks to Rui Fernandes and Stavros Papadakis for the discussions we had while I was writing this paper.

\section{PBG-groupoids}

Let us fix some conventions first. All manifolds considered 
in the paper are $C^{\infty}$ - differentiable, Hausdorff, paracompact and have a 
countable basis for their topology. We consider the arrows in a groupoid "from 
right to left'', i.e. the source $s\xi$ of an element $\xi$ in a groupoid 
$\Omega$ is considered to be on the right and its target $t\xi$ on the left. 
The object inclusion map of a groupoid $\Omega$ over a manifold $M$ is denoted 
by $1 : M \rightarrow \Omega$, namely $x \mapsto 1_{x}$. In \cite{LGLADG} 
the expression `Lie groupoid' stood for a locally trivial differentiable groupoid;
here we use the expression `transitive Lie groupoid' for clarity. 

\begin{df}\label{df:PBGgpd}
A {\em PBG-groupoid} is a Lie groupoid $\Upsilon \gpd P$ whose base is the total
space of a principal bundle $P(M,G)$ together with a right action of $G$ on 
the manifold $\Upsilon$ such that for all $(\xi,\eta) \in \Upsilon \times \Upsilon$ such that $s\xi = t\eta$ and $g \in G$ we have:
\begin{enumerate}
\item $t(\xi \cdot g) = t(\xi) \cdot g$ and $s(\xi \cdot g) =
s(\xi) \cdot g$
\item $1_{u \cdot g} = 1_{u} \cdot g$
\item $(\xi \eta) \cdot g = (\xi \cdot g)(\eta \cdot g)$
\item $(\xi \cdot g)^{-1} = \xi^{-1} \cdot g$
\end{enumerate}
\end{df}

We denote a PBG-groupoid $\Upsilon$ over the principal bundle $P(M,G)$ by 
$\Upsilon \gpd P(M,G)$ and the right-translation in $\Upsilon$ coming from the 
$G$-action by $\Tilde{R}_{g}$ for any $g \in G$. The right-translation in $P$ 
will be denoted by $R_{g}$. The properties in the previous definition show that 
the group $G$ of the base principal bundle $P(M,G)$ acts on $\Upsilon$ by 
automorphisms, namely $\Tilde{R}_{g}$ is an automorphism of the Lie groupoid 
$\Upsilon$ over the diffeomorphism $R_{g}$ for all $g \in G$. A morphism $\phi$ 
of Lie groupoids between two PBG-groupoids $\Upsilon$ and $\Upsilon'$ over the same 
principal bundle is called a morphism of PBG-groupoids, if 
$\phi \circ \Tilde{R}_{g} = \Tilde{R}_{g}' \circ \phi$ 
for all $g \in G$. In the same fashion, a {\em PBG-Lie group bundle} (PBG-LGB) is a 
Lie group bundle $F$ over the total space $P$ of a principal bundle $P(M,G)$ 
such that the group $G$ acts on $F$ by Lie group bundle automorphisms. We 
denote a PBG-LGB by $F \rightarrow P(M,G)$. It is easy to see that the gauge 
group bundle $I\Upsilon$ of a PBG-groupoid $\Upsilon \gpd P(M,G)$ is a PBG-LGB. 

Now let us describe the correspondence of transitive PBG-groupoids with extensions of Lie groupoids. For any given extension of Lie groupoids (\ref{gpdextn2}), choose a basepoint and take its corresponding extension of principal bundles 
\[
N \inj Q(M,H) \stackrel{\pi(id, \pi)}{\surj} P(M,G).
\]
That is to say that $\pi(id, \pi)$ is a surjective morphism of principal bundles and $H$ is an extension of the Lie group $G$ by $N$. 
This gives rise to the principal bundle $Q(P,N,\pi)$, which Mackenzie in \cite{Mackenzie:Cahiers} called the transverse bundle of the previous extension. The Lie groupoid $\Upsilon = \frac{Q \times Q}{N}$ corresponding to the transverse bundle admits the following action of $G$: 
\[
\langle v_2, v_1\rangle g = \langle v_2h, v_1h\rangle,
\]
for any $v_2, v_1 \in Q$ and $g \in G$, where $h\in H$ is any element with $\Tilde{\pi}(h) = g$. This action makes $\Upsilon \gpd P(M,G)$ a PBG-groupoid.
On the other hand, writing a transitive PBG-groupoid as an exact sequence $I\Upsilon\inj\Upsilon\surj P\times P$, 
the fact that $G$ acts by Lie groupoid automorphisms allows us to quotient the
sequence over $G$ and obtain the extension of Lie groupoids
\begin{equation}\label{eq:quot}
\frac{I\Upsilon}{G}\inj\frac{\Upsilon}{G}\surj \frac{P\times P}{G}. 
\end{equation}
The two processes are mutually inverse (for the details of the proof see \cite{Mackenzie:Cahiers}). This establishes the following result:
\begin{prop}\label{PBG:equiv:extn}
Any transverse PBG-groupoid $\Upsilon \gpd P(M,G)$ corresponds exactly to an extension of Lie groupoids $F \inj \Omega \surj \frac{P \times P}{G}$.
\end{prop}

Let us give now some examples of PBG-groupoids, which will also be useful in later sections.
\begin{ex}\label{ex:trPBGgpd}
{\em Consider a principal bundle $P(M,G)$ and a Lie group $H.$ Suppose given an 
action by automorphisms of $G$ on $H,$ say $(g,h) \mapsto R_{g}(h)$ for all 
$g \in G$ and $h \in H.$ That is to say that $R_{g} : H \rightarrow H$ is an 
automorphism of $H$ for all $g \in G.$ Form the trivial groupoid 
$P \times H \times P \gpd P(M,G).$ This is easily seen to be a transitive PBG 
groupoid. We will refer to it as the {\em trivial PBG-groupoid} 
corresponding to the given action.}
\end{ex}
\begin{ex}\label{ex:actPBGgpd}
{\em For this example, first let us recall the definition of an action groupoid. 
When a Lie group $G$ acts on a manifold $M$ the action groupoid $M \ract G$ is 
the product manifold $M \times G$ with the following groupoid structure: 
$s(x,g) = x,~t(x,g) = xg,~1_{x} = (x,e_{G})$. The multiplication is 
defined by $(xg,h) \cdot (x,g) = (x,gh)$ and the inverse of $(x,g)$ is 
$(xg,g^{-1})$. The action groupoid is transitive if and only if the action of 
$G$ on $M$ is transitive. Now consider a principal bundle $P(M,G)$. There is an 
action of $G$ on $P \ract G$ defined by $\Tilde{R}_{g}(u,h) = (ug,g^{-1}hg)$. 
This action makes the action groupoid a PBG-groupoid $P \ract G \gpd P(M,G)$.} 
\end{ex}
\begin{ex}\label{ex:frPBGgpd}
{\em A PBG vector bundle is a vector bundle $E$ over the total space $P$ of a 
principal bundle $P(M,G)$ such that the structure group $G$ acts on $E$ by vector 
bundle automorphisms. We denote a PBG vector bundle by $E \rightarrow P(M,G)$. 
The isomorphisms between the fibers of an arbitrary vector bundle define a 
transitive Lie groupoid called the frame groupoid. This is the groupoid 
corresponding to the frame principal bundle induced by the vector bundle under 
consideration. The frame groupoid $\Phi(E)$ of a PBG vector bundle 
$E \rightarrow P(M,G)$ has a canonical PBG structure over $P(M,G)$. Namely, for 
any $g \in G$, an isomorphism $\xi : E_{u} \rightarrow E_{v}$ between two fibers  
defines an isomorphism $\Tilde{R}_{g}(\xi) : E_{ug} \rightarrow E_{vg}$ by
\[
\Tilde{R}_{g}(\xi)(V) = \xi(V \cdot g^{-1}) \cdot g
\]
for all $V \in E_{ug}$. Moreover, there is a canonical morphism of PBG-groupoids 
$\epsilon : P \ract G \rightarrow \Phi(E)$ over the principal bundle $P(M,G)$. 
Namely, $\epsilon(u,g)$ is the isomorphism 
$\Hat{R}_{g} : E_{u} \rightarrow E_{ug}$ induced by the action of $G$ on $E$. 
It is straightforward to show that $\epsilon$ is indeed a morphism of Lie 
groupoids. It also preserves the actions because:
\begin{multline*}
\epsilon((u,h) \cdot g) = \epsilon(u \cdot g,g^{-1}hg) = 
[(\Hat{R}_{g^{-1}hg})_{u \cdot g} : E_{u \cdot g} \rightarrow E_{u \cdot hg}] = \\
= \Tilde{R}_{g}[(\Hat{R}_{h})_{u} : E_{u} \rightarrow E_{u \cdot h}] = 
\Tilde{R}_{g}(\epsilon(u,h)).
\end{multline*}}
\end{ex}
Let us now give an example of a PBG-groupoid that arises from an extension of 
principal bundles.
\begin{sloppypar}
\begin{ex}
{\em We consider the principal bundles
$SU(2)(S^{2},U(1),p)$ and $SO(3)(S^{2},SO(2),p')$. For the first bundle, 
denote a typical element
$\smallskip \begin{bmatrix}
s & t \\
-\Bar{t} & \Bar{s} \\
\end{bmatrix}$
of $SU(2)$ such that $|s|^{2} + |t|^{2} = 1$ by 
$\smallskip (s,t)$. Regard $U(1)$ as a subgroup of $SU(2)$ by mapping 
every $z \in U(1)$ to $(z,0)$ and let $p$ be
\[
(s,t) \mapsto 
(- 2Re(s \cdot t),- 2Im(s \cdot t),1 - 2 \cdot|t|^{2}).
\]
For the second bundle, regard $SO(2)$ as a subgroup of $SO(3)$ by $A \mapsto
\begin{bmatrix}
A & 0 \\
0 & 1 \\
\end{bmatrix}$
and let $\smallskip p'$ be $A \mapsto A \cdot e_{3}$ where $(e_{1},e_{2},e_{3})$ 
is the usual basis of $\reals^{3}$. We define a morphism between these  principal 
bundles making use of quaternions. Consider any 
$q = (s, t) \in SU(2)$ as a unit quaternion and every
$r = (r_{1},r_{2},r_{3}) \in \reals^{3}$ as a vector quaternion. 
Then $q \cdot r \cdot q^{-1}$ is a vector quaternion, i.e. 
$q \cdot r \cdot q^{-1} \in \reals^{3}.$ For every
$q = (s,t) \in SU(2)$ we define 
$A_{q} \co \reals^{3} \rightarrow \reals^{3}$ as
$A_{q}(r) = q \cdot r \cdot q^{-1}$. This is an element in $SO(3)$. Let 
$\phi : U(1) \rightarrow SO(2)$ be the restriction of $A$ to $U(1)$. If we write 
a $z \in U(1)$ as $z = e^{\iota \theta} = \cos{\theta} + \iota \sin{\theta}$. for 
some $\theta \in \reals$, then it is easy to see that
\[
\phi(z) =
\begin{bmatrix}
\cos(2\theta) & \sin(2\theta) \\
- \sin(2\theta) & \cos(2\theta)
\end{bmatrix}
\in SO(2).
\]
This shows both that $\phi$ is a surjective submersion and that its kernel 
is $\integers_{2}$. We therefore get the extension 
\[
\integers_{2} \inj SU(2)(S^{2},U(1))\stackrel{R(id_{S^{2}},\phi)}{\surj} 
SO(3)(S^{2},SO(2)).
\]
Its transverse bundle is
$SU(2)(SO(3),\integers_{2},A)$ and the PBG-groupoid it induces is 
\[
\frac{SU(2) \times SU(2)}{\integers_{2}} \gpd SO(3)(S^{2},SO(2)),
\]
where the action of $SO(2)$ on $\frac{SU(2) \times SU(2)}{\integers_{2}}$ is:
\[
\langle (s_{1},t_{1}),(s_{2},t_{2})\rangle \cdot z = 
\langle (s_{1}z,t_{1}\Bar{z}),(s_{2}z,t_{2}\Bar{z}) \rangle.
\]
It was shown in \cite[II\S7]{LGLADG} that for any Lie group $G$ and closed 
subgroup $H \leq G$ the gauge groupoid $\frac{G \times G}{H}$ is isomorphic to 
the action groupoid $\frac{G}{H}\ract G,$ where $G$ acts on $\frac{G}{H}$ 
via $(g,g'H) \mapsto (g \cdot g')H$. The isomorphism is 
\[
\langle g_{1},g_{2} \rangle  \mapsto (g_{1} \cdot g_{2}^{-1},g_{2}H). 
\]
\begin{sloppypar}
Therefore $\smallskip \frac{SU(2) \times SU(2)}{\integers_{2}} \gpd 
SO(3)(S^{2},SO(2))$ is isomorphic to the action PBG-groupoid 
$\smallskip \frac{SU(2)}{\integers_{2}}\ract SU(2) \gpd SO(3)(S^{2},SO(2))$. 
Note that $\smallskip \frac{SU(2)}{\integers_{2}} \cong \reals P^{3}$. 
The action of $SO(2)$ on $\reals P^{3} \ract SU(2)$ is 
\end{sloppypar}
\[
((w_{1},w_{2}),(s , t) \cdot \integers_{2}) \cdot R_{\theta} =
((w_{1},w_{2}),(s \cdot z,t \cdot \Bar{z}) \cdot \integers_{2}).
\]}
\end{ex}
\end{sloppypar}

\section{PBG-algebroids}

This is an account of PBG structures on the algebroid level.
We show that every PBG-Lie algebroid corresponds to an extension of Lie algebroids; the account in \cite{Mackenzie:Cahiers} gave only a partial result of this type.  Some fundamental results are also proved, which will be useful in the study of the connection theory in later sections of this paper. Let us begin by describing the notion of a PBG structure on the algebroid 
level. We start with the special case of Lie algebra bundles.
\begin{df}\label{df:PBGLAB}
A {\em PBG-Lie algebra bundle} (PBG-LAB) is a Lie algebra bundle $L$ over the 
total space $P$ of a principal bundle $P(M,G)$, together with an action of $G$ on 
$L$ such that each right-translation $\Bar{R}_{g} : L \rightarrow L$ is a Lie 
algebra bundle automorphism over the right translation $R_{g} : P \rightarrow P$.
\end{df}
We denote a PBG-LAB by $L \rightarrow P(M,G)$. A morphism between two PBG-LABs 
over the same principal bundle is a morphism of LABs which preserves the group 
actions. 
\begin{ex}\label{ex:trPBGLAB}
{\em Let $P(M,G)$ be a principal bundle and $\hoh$ a Lie algebra which admits a 
right action $(V,g) \mapsto V \cdot g$ of $G$ such that 
\[
[V,W] \cdot g = [V \cdot g,W \cdot g]
\]
for all $V,W \in \hoh$. Now the trivial LAB $P \times \hoh$ admits the 
$G$-action $\Bar{R}_{g}(u,V) = (ug,V \cdot g)$. This makes it a PBG-LAB, the 
{\em trivial PBG-LAB}.}
\end{ex}
\begin{df}
A {\em Lie algebroid} is a vector bundle $A$ on base $M$ together with a vector bundle map $q : A \rightarrow TM$, called the {\em anchor} of $A$, and a bracket 
$[\ ,\ ] \co \Gamma A \times \Gamma A \rightarrow \Gamma A$ which is $\reals$-bilinear, alternating, satisfies the Jacobi identity, and is such that
\begin{enumerate}
\item $q([X,Y]) = [q(X),q(Y)]$,
\item $[X,uY] = u[X,Y] + q(X)(u)Y$
\end{enumerate}
for all $X,Y \in \Gamma A$ and $u \in C^{\infty}(M)$. 
\end{df}

Given a trivial PBG-LAB $P \times \hoh \rightarrow P(M,G)$, recall that the 
Whitney sum vector bundle $TP \oplus (P \times \hoh)$ has a trivial Lie algebroid 
structure over $P$. The anchor map is the first projection and the Lie bracket 
on its sections is defined by the formula
\[
[X \oplus V,Y \oplus W] = [X,Y] \oplus \{ X(W) - Y(V) + [V,W] \}
\]
for all $X,Y \in \Gamma TP$ and smooth $\hoh$-valued functions on $V,W$ on $P$. 
The $G$-action on $P \times \hoh$ defines a PBG structure on this Lie algebroid. 
Namely, with the notation used in the previous example (\ref{ex:trPBGLAB}) define 
the following action of $G$ on $TP \oplus (P \times \hoh)$:
\[
(X \oplus (u,V),g) \mapsto \Hat{R}_{g}(X \oplus (u,V)) = T_{u}R_{g}(X) \oplus (ug,V \cdot g).
\]
Denote by $\Hat{R}^{\Gamma}_{g} \co \Gamma(TP \oplus (P \times \hoh)) \rightarrow 
\Gamma(TP \oplus (P \times \hoh))$ the corresponding action on the sections of 
$TP \oplus (P \times \hoh)$. This is given by the formula
\[
\Hat{R}^{\Gamma}_{g}(X \oplus V)_{u} = 
T_{ug^{-1}}R_{g}(X_{ug^{-1}}) \oplus (V_{ug^{-1}} \cdot g)
\]
for all $X \in \Gamma TP,~V \in C^{\infty}(P,\hoh)$ and $u \in P$. This action 
preserves the Lie bracket in $\Gamma(TP \oplus (P \times \hoh))$ and is an action 
by automorphisms on the trivial Lie algebroid. This construction is an example
of the notion of a PBG-algebroid.

\begin{df}\label{df:PBGalgd}
A {\em PBG-algebroid} over the principal bundle $P(M,G)$ is a Lie 
algebroid $A$ over $P$ together with a right action of $G$ on $A$ denoted by 
$(X,g) \mapsto \Hat{R}_{g}(X)$ for all $X \in A,~ g \in G$ such that each 
$\Hat{R}_{g} : A \rightarrow A$ is a Lie algebroid automorphism over the right 
translation $R_{g}$ in $P$. 
\end{df}

We denote a PBG-algebroid $A$ over $P(M,G)$ by $A \Rightarrow P(M,G)$. The action 
of $G$ on $A$ induces an action of $G$ on $\Gamma A$, namely
\[
X \cdot g = \Hat{R}_{g} \circ X \circ R_{g^{-1}}.
\]
for all $g \in G$ and $X \in \Gamma A$. The right-translation with respect to 
this action is denoted by $\Hat{R}_{g}^{\Gamma} \co \Gamma A \rightarrow \Gamma A$ 
for all $g \in G$. With this notation definition \ref{df:PBGalgd} implies that 
\[
\Hat{R}_{g}^{\Gamma}([X,Y]) = [\Hat{R}^{\Gamma}_{g}(X),\Hat{R}^{\Gamma}_{g}(Y)] 
\]
for all $X,Y \in \Gamma A$ and $g \in G$. The following result strengthens 3.2 of \cite{Mackenzie:Cahiers}, where the quotient manifold was assumed to exist. 

\begin{prop}
\label{thm:strong}
Let $A$ be a transitive PBG-Lie algebroid on $P(M,G)$. Then the quotient manifold
$A/G$ exists and inherits a quotient structure of transitive Lie algebroid
from $A$; further, it is an extension
$$
\frac{L}{G}\inj \frac{A}{G}\surj \frac{TP}{G} 
$$
of the Lie algebroid of the gauge groupoid of $P(M,G)$ by the quotient
LAB $L/G$. 
\end{prop}

\pf
The main requirement is to prove that the manifold $A/G$ exists. We apply
the criterion of Godement \cite{Dieudonne}. Denote the projection $A\to P$ 
by $p_A$, and write
$$
\Ga' = \{(X, Xg)\st X\in A,\ g\in G\};
$$
we must show that $\Ga'$ is a closed submanifold of $A\times A$. Now
$\Ga'\subseteq (p_A\times p_A)^{-1}(\Ga)$ where 
$\Ga = \{(u, ug)\st u\in P,\ g\in G\}.$ Since $\Ga$ is a closed submanifold
of $P\times P$, and $p_A\times p_A$ is a surjective submersion, it suffices to
prove that $\Ga'$ is a closed submanifold of $(p_A\times p_A)^{-1}(\Ga)$. Define
$$
f\co (p_A\times p_A)^{-1}(\Ga)\to A,\qquad (X, Y)\mapsto Xg - Y
$$
where $p_A(Y) = p_A(X)g$. From the local triviality of $A$, it easily follows
that $f$ is a surjective submersion. The preimage of the zero section under
$p_A$ is $\Ga'$, and this shows that $\Ga'$ is a closed submanifold. Denote the
quotient projection $A\to A/G$ by $\natural$. 

The vector bundle structure of $A$ quotients to $A/G$ in a straightforward
fashion. Since the anchor $q\co A\to TP$ is $G$--equivariant, it quotients to
a vector bundle morphism $A/G\to TP/G$ which is again a surjective submersion; 
denote this by $\pi$, and define $r = \Tilde{q}\circ \pi$ where $\Tilde{q}$ is 
the anchor of $TP/G$. 

For the bracket structure of $A/G$, note first that $\Ga(A/G)$ can be 
identified with the $C^\infty(M)$ module of $G$--equivariant sections of $A$
as in the case of the Atiyah sequence of a principal bundle. 
Since the bracket on $\Ga A$ restricts to the $G$--equivariant sections by
assumption, this bracket transfers to $\Ga(A/G)$. It is now straightforward
to check that this makes $A/G$ a Lie algebroid on $M$ with anchor $r$, and
$\natural\co A\to A/G$ a Lie algebroid morphism over $p$. That $\pi$ is a
Lie algebroid morphism with kernel $L/G$ is easily checked. 
\boom

It is a straightforward exercise to verify that PBG-groupoids differentiate to PBG-algebroids. 
The next result, which appears in \cite{Mackenzie:Cahiers}, gives a converse. 

\begin{thm}
Let $\Ups$ be an $s$--simply connected locally trivial Lie groupoid
with base the total space of a principal bundle $P(M, G)$. Suppose that 
for all $g\in G$ there is given a Lie algebroid automorphism 
$\Tilde{R}_g\co A\Ups\to A\Ups$ which defines the structure of a 
PBG-algebroid on $A\Ups$. Then there is a natural structure of 
PBG-groupoid on $\Ups$ which induces on $A\Ups$ the given PBG-Lie
algebroid structure. 
\end{thm}

The examples which follow will be useful later on.

\begin{ex}\label{ex:notsimplyconnected}
{\em Consider a transitive Lie algebroid $A$ over a manifold $M$ such that 
$\pi_{1}(M) \neq 0$. Let $\Tilde{M}$ be the covering space of $M$ and 
$p \co \Tilde{M} \rightarrow M$ the covering projection. Then we have the principal 
bundle $\Tilde{M}(M,\pi_{1}(M),p)$. Denote $p^{!!}A$ the pullback
$$
\begin{diagram}
p^{!!}A & \rto^{p^*} & p^*A \\
\dto<{q'} &  & \dto>{p^*(q)} \\
T\Tilde{M} & \rto^{T(p^{*})} & p^*TM \\
\end{diagram}
$$
Sections of $p^{!!}A$ are of the form $X' \oplus C$, where $X' \in \Gamma T\Tilde{M}$, $C \in \Gamma (p^*A)$ and
\[
T(p^*)(X') = p^*(q)(C).
\]
This construction is described in \cite{Mackenzie:constronalgds}, but let us recall here the Lie algebroid structure of $p^{!!}A$. If we write $C = \Sigma u_{i}(X_i \circ p)$ with $u_i \in C(\Tilde{M})$ and $X_{i} \in \Gamma A$, then 
\[
p^*(q)(C) = \Sigma u_{i}(q(X_i \circ p))
\]
with $T(p) \circ X' = \Sigma u_i(q(X_i \circ p))$. We define a 
Lie algebroid structure on $p^{!!}A$ with $q'$ as the anchor map and bracket given by
\begin{multline*}
[X'\oplus \Sigma u_i(X_i \circ p), Y' \oplus \Sigma u_j(Y_j \circ p)] = \\
= [X',Y'] \oplus \{ \Sigma u_i u_j ([X_i,Y_j] \circ p) + \Sigma X'(u_j)(Y_j \circ p) 
- \Sigma Y'(u_i)(X_i \circ p) \}.
\end{multline*}
Now the pullback Lie algebroid $p^{!!}A$ is a PBG-algebroid over $\Tilde{M}(M,\pi_{1}(M),p)$.} 
\end{ex}
\begin{ex}\label{ex:actPBGalgd}
{\em Consider the action PBG-groupoid $P \ract G \gpd P(M,G)$ constructed in 
\ref{ex:actPBGgpd}. This differentiates to the {\em action PBG-algebroid} 
$P \ract \gog \Rightarrow P(M,G)$. Here $P \ract \gog$ is the product manifold 
$P \times \gog$ and the anchor is the map that associates to every $V \in \gog$ 
the fundamental vector field $V^{\dagger} \in \Gamma TP$. The Lie bracket is 
given by the formula
\begin{eqnarray}\label{eqn:ptwsbracket}
[V,W] = V^{\dagger}(W) - W^{\dagger}(V) + [V,W]^{\ptwise},
\end{eqnarray}
where $[V,W]^{\ptwise}$ stands for the point-wise bracket in $\gog$. The 
$G$-action that makes it a PBG-algebroid is of course
\[
\Hat{R}_{g}(u,V) = (ug,\Ad_{g^{-1}}(V)).
\]}
\end{ex}
\begin{ex}\label{ex:CDO[L]}
{\em Given a PBG-LAB $L \rightarrow P(M,G)$, let $\Phi[L]$ denote the groupoid 
of Lie algebra isomorphisms between the fibers of $L$. This is a PBG-groupoid over 
$P(M,G)$ in the same way as $\Phi(L)$ (see \ref{ex:frPBGgpd}). It was shown in 
\cite{LGLADG} that $\Phi[L]$ differentiates to the Lie algebroid $\CDO[L]$ over 
$P$. The notation $\CDO$ stands for ``covariant differential operator''. It is a transitive Lie algebroid and its 
sections are those first or zeroth order differential operators 
$D \co \Gamma L \rightarrow \Gamma L$ such that:
\begin{enumerate}
\item For all $D \in \Gamma\CDO[L]$ there is a vector field 
$\sharp(D) \in \Gamma TP$ such that
\[
D(f\mu) = fD(\mu) + \sharp(D)(f)\mu
\]
for every $\mu \in \Gamma L$ and $f \in C^{\infty}(P)$.
\item The operators $D$ act as derivations of the bracket, i.e. 
\[
D([\mu_{1},\mu_{2}]) = [D(\mu_{1}),\mu_{2}] + [\mu_{1},D(\mu_{2})]
\]
for all $\mu_{1},\mu_{2} \in \Gamma L$. 
\end{enumerate}
The anchor map of $\CDO[L]$ is exactly the map $\sharp$ established in (i). The 
adjoint bundle of this Lie algebroid is the LAB of endomorphisms of $L$ which
are derivations of the bracket. We therefore have the exact sequence
\[
\Der(L) \inj \CDO[L] \stackrel{\sharp}{\surj} TP.
\]
The action of $G$ on $\Phi[L]$ differentiates to the action 
$(D,g) \mapsto R_{g}^{\CDO}(D)$ on the sections of $\CDO[L]$ defined by
\[
[R_{g}^{\CDO}(D)](\mu) = 
\Bar{R}_{g} \circ D(\Bar{R}_{g^{-1}} \circ \mu \circ R_{g}) \circ R_{g^{-1}}
\]
for all $g \in G, D \in \Gamma\CDO[L]$ and $\mu \in \Gamma L$. (Recall that 
$\Bar{R}_{g}$ denotes the right-translation on the sections of $L$.) This action 
makes $\CDO[L]$ a PBG-algebroid over $P(M,G)$.} 
\end{ex}

In general, $\CDO[L]$ is the analogue of the automorphism group in the Lie algebroid framework. In the PBG setting, we have the following definition:
\begin{df}\label{df:PBGalgdrep}
Let $A \Rightarrow P(M,G)$ be a (transitive) PBG-algebroid and $K \rightarrow P(M,G)$ a PBG-LAB. An {\em equivariant representation} of $A$ on $K$ is an equivariant Lie algebroid morphism $\rho\co A \rightarrow \CDO[L]$.
\end{df}

\begin{prop}\label{prop:adjbundle}
Let $A \Rightarrow P(M,G)$ be a transitive PBG-algebroid. Then its adjoint bundle 
$L$ is a PBG-LAB and $\CDO[L]$ a PBG-algebroid, both over the principal bundle 
$P(M,G)$.
\end{prop}

\pf~Since $A$ is transitive we have the exact sequence of Lie algebroids 
\[
L \stackrel{j}{\inj} A \stackrel{\sharp}{\surj} TP. 
\]
Let $(V,g) \mapsto \Bar{R}_{g}(V)$ denote the action of $G$ on $\Gamma L$ 
defined by $j(\Hat{R}^{\Gamma}_{g}(V)) = \Hat{R}_{g}^{\Gamma}(j(V))$. The result 
now follows from the fact that $A$ is a PBG-algebroid. \boom

A Lie algebroid morphism $\phi : A \rightarrow A'$ between two PBG-algebroids 
over the same principal bundle $P(M,G)$ is a {\em morphism of PBG-algebroids} if 
$\phi \circ \Hat{R}_{g} = \Hat{R}_{g} \circ \phi$ for all $g \in G$. It is easy 
to verify that morphisms of PBG-groupoids differentiate to PBG-algebroids. 
\begin{prop}\label{prop:integrPBGmorph}
Let $\Upsilon$ and $\Upsilon'$ be PBG-groupoids over the same principal bundle 
$P(M,G)$ such that $\Upsilon$ is $s$-simply connected. Then every morphism of 
PBG-algebroids $\phi_{*} : A\Upsilon \rightarrow A\Upsilon'$ integrates to a unique 
morphism of PBG-groupoids $\phi \co \Upsilon \rightarrow \Upsilon'$.
\end{prop}
\pf~It was shown in \cite{Mackenzie:Liebialgds} that the Lie algebroid morphism 
$\phi_{*}$ integrates uniquely to a morphism of Lie groupoids 
$\phi \co \Upsilon \rightarrow \Upsilon'$. It suffices to show that $\phi$ respects 
the $G$-actions. For every $g \in G$ let $\Tilde{R}_{g}$ and $\Tilde{R}_{g}'$ be 
the right-translations induced by the actions of $G$ on $\Upsilon$ and $\Upsilon'$ 
respectively. These differentiate to the right-translations $\Hat{R}_{g}$ and 
$\Hat{R}_{g}'$ respectively on the Lie algebroid level. Since $\phi_{*}$ is a 
morphism of PBG-algebroids, for every $g \in G$ we have 
$\Hat{R}_{g}' \circ \phi \circ \Hat{R}_{g^{-1}} = \phi_{*}$. The morphism of Lie 
groupoids $\Hat{R}_{g}' \circ \phi \circ \Hat{R}_{g^{-1}}$ differentiates to 
$\Hat{R}_{g}' \circ \phi_{*} \circ \Hat{R}_{g}$. The uniqueness of $\phi$ yields 
$\Tilde{R}_{g}' \circ \phi \circ \Tilde{R}_{g^{-1}} = \phi$ for all 
$g \in G$. \boom

\begin{ex}\label{ex:derreps}
{\em The natural morphism of PBG-groupoids arising from a PBG-LAB 
$L \rightarrow P(M,G)$ differentiates to a morphism of PBG-algebroids 
$\epsilon_{*} : P \ract \gog \rightarrow \CDO[L]$. Such morphisms are known 
as ``derivative representations'' in the context of Poisson geometry. We will 
refer to this one as the {\em natural derivative representation} induced by $L$.}
\end{ex}

A Lie subalgebroid $A'$ of a PBG-algebroid $A$ is a PBG-subalgebroid if the 
inclusion $A' \inj A$ is a morphism of PBG-algebroids. 
\begin{prop}\label{prop:adjPBGalgd}
If $A\Rightarrow P(M,G)$ is a transitive PBG-algebroid with adjoint bundle 
$L$ then $\ad(A)$ is a PBG-subalgebroid of $\CDO[L]$.
\end{prop}

\pf~The Lie algebroid $\ad(A)$ was shown to be a Lie subalgebroid of 
$\CDO[L]$ in \cite{LGLADG}. It suffices to show that it is closed under 
the $G$-action defined in \ref{ex:CDO[L]}. Indeed, suppose $X \in \Gamma A$. 
Then:
\begin{multline*}
[R_{g}^{\CDO}(\ad_{X})](V) = 
\Bar{R}_{g} \circ \ad_{X}(\Bar{R}_{g^{-1}} \circ V \circ R_{g}) \circ R_{g^{-1}} 
= \\
= \Bar{R}_{g} \circ [X,\Bar{R}_{g^{-1}} \circ V \circ R_{g}] \circ R_{g^{-1}} = 
[\Bar{R}_{g}(X),V] = \ad_{\Bar{R}_{g}(X)}(V)
\end{multline*}
for all $V \in \Gamma L$. So, $R_{g}^{\CDO}(\ad_{X}) \in \ad(A)$. \boom

The Lie algebroid $\ad(L)$ is integrable as a subalgebroid of the integrable 
Lie algebroid $\CDO[L]$. We denote the $\alpha$-connected Lie groupoid it 
integrates to by $\Int(A)$, following the notation of \cite{LGLADG}. The proof 
of the following statement is immediate.

\begin{cor}\label{cor:Int{A}}
If $A \Rightarrow P(M,G)$ is a transitive PBG-algebroid then
 $\Int{A} \gpd P(M,G)$ is a PBG-subalgebroid of $\Phi[L]$. 
\end{cor}

We can also consider $\ad(A)$ as the image of the adjoint representation 
$\ad : A \rightarrow \CDO[L]$. In this sense, $\ad$ is a representation of PBG 
algebroids, i.e. a morphism of PBG-algebroids. 
\begin{df}\label{df:PBGalgdreps}
Let $A$ be a PBG-algebroid and $L$ a PBG-LAB, both over the same principal 
bundle $P(M,G)$. A representation $\rho : A \rightarrow \CDO[L]$ of $A$ in $L$ is 
a {\em representation of PBG-algebroids} if it is a morphism of PBG-algebroids.
\end{df}
In particular, for all $g \in G, X \in \Gamma A$ and $\mu \in \Gamma L$ it is 
required:
\[
[R_{g}^{\CDO}(\rho(X))](\mu) = \rho(\Hat{R}_{g}^{\Gamma}(X))(\mu).
\]
The developed form of this formula is: 
\[
\Bar{R}_{g} \circ \rho(X)(\Bar{R}_{g^{-1}} \circ \mu \circ R_{g}) 
= \rho(\Hat{R}_{g}^{\Gamma}(X))(\mu) \circ R_{g}.
\]

\section{Infinitesimal connections for PBG structures}

Lie groupoids and Lie algebroids provide a natural framework for the study of 
connection theory. For example, given a principal bundle $P(M,G,p)$, consider its 
corresponding Atiyah sequence, namely the exact sequence of vector bundles
\[
\frac{P \times \gog}{G} \inj \frac{TP}{G} \stackrel{p_{*}}{\surj} TM
\]
(where the action of $G$ on $\gog$ implied is simply the adjoint), 
together with the bracket of sections of the vector bundle $\frac{TP}{G}$
which is obtained by identifying them with $G$--invariant vector fields on $P$. 
This bracket is preserved by the vector bundle morphism $p_{*}$. Thus the Atiyah 
sequence of a principal bundle is a transitive Lie algebroid. In this setting, infinitesimal connections of the principal bundle $P(M,G,p)$ 
correspond exactly to the right-splittings of the Atiyah sequence. The 
holonomy of such connections is studied using the global object corresponding 
to the Atiyah sequence and that is the transitive Lie groupoid corresponding 
to the principal bundle we started with, namely the quotient manifold 
$\frac{P \times P}{G}$ (over the diagonal action). 

Consider a transitive PBG-algebroid
\begin{eqnarray}\label{PBGalgdextn}
L\inj A \surj TP
\end{eqnarray}
over the principal bundle $P(M,G)$ and form the associated extension of Lie algebroids
\begin{eqnarray}\label{assoclalgdextn}
\frac{L}{G}\inj \frac{A}{G} \surj \frac{TP}{G}.
\end{eqnarray}
We will write $\Phi = \frac{P\times P}{G}$ and $B = \frac{A}{G},\
K = \frac{L}{G}$. As vector bundles $L \isom p^*K,\ A \isom p^*B$ and
$TP \isom p^*(A\Phi)$. These are actually isomorphisms of Lie algebroids
with respect to the action Lie algebroid structures on the pullback bundles;
see \cite{Mackenzie:Cahiers} for this. A connection (or an infinitesimal connection if emphasis
is required) in $A$ is a right--inverse $\gamma\co TP\to A$ to (\ref{PBGalgdextn}). If
$\gamma$ is equivariant with respect to the $G$ actions, that is, if 
\begin{equation}
\label{eq:isomet}
\gamma \circ TR_{g} = \Hat{R}_{g} \circ \gamma
\end{equation}
for all $g\in G$, then $\gamma$ quotients to a linear map $\gamma^G\co A\Phi\to B$
which is right--inverse to (\ref{assoclalgdextn}). Conversely, given a right--inverse $\sigma
\co A\Phi\to B$ to (\ref{assoclalgdextn}), the pullback map $p^*\sigma$ is a connection in $A$ 
which satisfies (\ref{eq:isomet}). Connections in $A$ satisfying (\ref{eq:isomet})
are not equivariant in the standard sense of, for example, \cite{GHV}, and 
in order to distinguish from the standard notion we call these 
connections {\em isometablic}. This is an adaptation of the Greek translation 
of the Latin word ``equivariant''. 

\begin{df}\label{conns:PBG}
Let $A \Rightarrow P(M,G,p)$ be a transitive PBG-algebroid. A connection 
$\gamma : TP \rightarrow A$ is called {\em isometablic}, if it satisfies 
(\ref{eq:isomet}). 
\end{df}

Writing a transitive PBG-algebroid as an extension of PBG-algebroids $L \stackrel{\iota}{\inj} A \stackrel{\sharp}{\surj} TP$, the curvature of any connection $\gamma$ in $A$ is the 2-form $\Omega_{\gamma} : TP \times TP \rightarrow L$ defined by 
\[
\Omega_{\gamma}(X,Y) = \gamma([X,Y]) - [\gamma(X),\gamma(Y)].
\]
If $\gamma$ is isometablic then its curvature also preserves the actions, namely for all $g \in G$ we have
\[
\Omega_{\gamma} \circ (TR_{g} \times TR_{g}) = \Bar{R}_{g} \circ \Omega_{\gamma}.
\]
An isometablic {\em back connection} of $A \Rightarrow P(M,G)$ is a morphism of vector bundles $\omega : A \rightarrow L$ such that $\omega \circ \iota = id_{L}$ 
and $\iota \circ \omega \circ \Hat{R}_{g} = \Bar{R}_{g} \circ \iota \circ \omega$ for all $g \in G$. Isometablic connections correspond to isometablic back 
connections via the formula
\[
\iota \circ \omega + \gamma \circ \sharp = id_{A}.
\]
\begin{ex}\label{ex:stflconn}
Consider the PBG-algebroid $TP \oplus (P \times \hoh)$ over a principal bundle $P(M,G)$ constructed in Section 3. The connection 
$\gamma^{0} : TP \rightarrow TP \oplus (P \times \hoh)$ defined by $X \mapsto X \oplus 0$ for all $X \in TP$ is isometablic and flat. This is  the {\em standard flat connection}.
\end{ex}\label{ex:flatisomconn}
If $\phi : A \rightarrow A'$ is a morphism of PBG-algebroids over the same principal bundle and $\gamma$ is an isometablic connection in $A$ then 
$\phi \circ \gamma$ is also an isometablic connection in $A'$.
Any connection in a transitive PBG-algebroid $A \Rightarrow P(M,G)$ gives rise to a Koszul connection $\nabla$ of the adjoint bundle $L$. Namely, define 
$\nabla^{\gamma} \co \Gamma TP \times \Gamma L \rightarrow \Gamma L$ by 
\[
\nabla^{\gamma}_{X}(V) = [\gamma(X),\iota(V)]
\]
for all $X \in \Gamma TP$ and $V \in \Gamma L$. This is the {\em adjoint connection} induced by $\gamma$. If $\gamma$ is an isometablic connection, then its 
adjoint connection satisfies
\begin{multline*}
\Bar{R}_{g}(\nabla^{\gamma}_{X}(V)) = \Bar{R}_{g}([\gamma(X),\iota(V)]) = [\Bar{R}_{g}(\gamma(X)),\Bar{R}_{g}(\iota(V))] = \\
= [\gamma(TR_{g}(X)),\iota(\Bar{R}_{g}(V))] = \nabla^{\gamma}_{TR_{g}(X)}(\Bar{R}_{g}(V))
\end{multline*}
for all $X \in \Gamma TP$ and $V \in \Gamma L$. 
The curvature of $\gamma$ satisfies:
\[
\odot \{ \nabla^{\gamma}_{X}(\Omega_{\gamma}(Y,Z)) - \Omega_{\gamma}([X,Y],Z) \} = 0
\]
for all $X,Y,Z \in \Gamma TP$. Here $\odot$ denotes the sum over all the permutations of $X,Y$ and $Z$. This is the (second) {\em Bianchi identity}.

\begin{prop}\label{basisom}
Suppose given a transitive PBG-algebroid $A \Rightarrow P(M,G,p)$ and consider its corresponding extension of Lie algebroids over $M$. The connections of the (transitive) Lie algebroid $\frac{A}{G} \rightarrow M$ are equivalent to the isometablic connections of $A$ which vanish on the kernel $T^{p}P$ of $Tp : TP \rightarrow TM$.
\end{prop}

\pf Consider an isometablic connection $\gamma : TP \rightarrow A$ such that $\gamma(X) = 0$ if $X \in T^{p}P$. This quotients to a splitting $\gamma^{/G} : \frac{TP}{G} \rightarrow \frac{A}{G}$. Given a connection $\delta : TM \rightarrow \frac{TP}{G}$ of the principal bundle $P(M,G)$, define
\[
\Tilde{\gamma} = \gamma^{/G} \circ \delta : TM \rightarrow \frac{A}{G}.
\]
The assumption that $\gamma$ vanishes on the kernel of $Tp$ makes the definition of $\Tilde{\gamma}$ independent from the choice of $\delta$. It follows immediately from the assumption that $\delta$ is a connection of $P(M,G)$ and $\gamma^{/G}$ is a splitting of (\ref{assoclalgdextn}) that this is a connection of the Lie algebroid $\frac{A}{G}$.

Conversely, given a connection $\theta : TM \rightarrow \frac{A}{G}$ of the Lie algebroid $\frac{A}{G}$, compose it with the anchor map $p^* : \frac{TP}{G} \surj TM$ of the Atiyah sequence corresponding to the bundle $P(M,G,p)$ to the vector bundle morphism
\[
\Bar{\theta} = \theta \circ p^* : \frac{TP}{G} \rightarrow \frac{A}{G}.
\]
Denote $\natural : TP \rightarrow \frac{TP}{G}$ and $\natural^A : A \rightarrow \frac{A}{G}$ the natural projections. Since $\natural^A$ is a pullback over $p \co P \rightarrow M$, there is a unique vector bundle morphism $\gamma : TP \rightarrow A$ such that
\[
\natural^A \circ \gamma = \Bar{\theta} \circ \natural.
\]
Due to the $G$-invariance of $\natural$ and $\natural^A$the morphism of vector bundles $\Hat{R}_{g^{-1}} \circ \gamma \circ TR_{g}$ also satisfies the previous equation for every $g \in G$, therefore it follows from the uniqueness argument that $\gamma$ is isometablic. It is an immediate consequence of the previous equation that $\gamma$ vanishes at $T^pP$. 

To see that it is indeed a connection of $A$, let us recall the fact that $\theta$ is a connection of $\frac{A}{G}$. This gives $p^* \circ \sharp^{/G} \circ \theta = \id_{TM}$. Now $\sharp^{/G} = \natural \circ \sharp$ and by definition we have $p^{*} \circ \natural = Tp$, therefore $Tp \circ \sharp \circ \theta = \id_{TM}$. Now take an element $X \in TP$. Then $Tp(X) \in TM$, and it follows from this equation that there exists an element $g \in G$ such that 
\[
(\sharp \circ \theta)(Tp(X)) = X \cdot g.
\]
Multiplying this by $g^{-1}$ and using the $G$-invariance of $Tp$ we get
\[
\sharp \circ (\theta \circ Tp) = \id_{TP}.
\]
Finally, from the properties of the pullback, it follows immediately that $\gamma$ is the map $(\pi, \Bar{\theta} \circ \natural)$, where $\pi \co TP \rightarrow P$ is the natural projection of the tangent bundle. It is straightforward to check that this reformulates to $(\pi, \theta \circ Tp)$, and this proves that $\gamma$ is a connection.
\boom
\begin{cor}\label{locflatisomconn}
Let $A \Rightarrow P(M,G)$ be a transitive PBG-algebroid. A flat connection of the Lie algebroid $\frac{A}{G} \rightarrow M$ gives rise to a unique flat connection of $A$ which vanishes on the kernel of $Tp$.
\end{cor}

\begin{df}\label{df:isomKoszconn}
An {\em isometablic Koszul connection} of the PBG-LAB $L \Rightarrow P(M,G)$ is a vector bundle morphism $\nabla : TP \rightarrow \CDO[L]$ such that
\[
\nabla_{TR_{g}(X)}(\Bar{R}_{g}(V)) = \Bar{R}_{g}(\nabla_{X}(V))
\]
for all $g \in G,~X \in \Gamma TP$ and $V \in \Gamma L.$
\end{df}
Now let us give a few examples of such connections.
\begin{ex}\label{prop:stndflatisomconn}
Suppose $P \times \hoh \Rightarrow P(M,G)$ is a trivial PBG-LAB (in the sense of \ref{ex:trPBGLAB}).
Then the {\em Koszul connection} $\nabla^{0} \co \Gamma TP \rightarrow \Gamma \CDO[P \times \hoh]$
defined by $\nabla^{0}_{X}(V) = X(V)$ for every $X \in \Gamma TP,~V \in C^{\infty}(P,\hoh)$ is
isometablic and flat. To verify this, note that it is the adjoint connection of the standard flat connection $\gamma^0$ introduced in \ref{ex:stflconn}
\end{ex}
Further examples arise when we consider the $\Hom$ functor. Now we present some of them, which are necessary for the proof of the Ambrose-Singer theorem in Section 7. 

Suppose $L \rightarrow P(M,G)$ and $L' \rightarrow P(M,G)$ are PBG-Lie algebra bundles. The vector 
bundle $\Hom^{n}(L,L') \rightarrow P$ is a Lie algebra bundle with Lie bracket
\[
[\phi_{1},\phi_{2}](\mu_{1},\ldots,\mu_{n}) = [\phi_{1}(\mu_{1},\ldots,\mu_{n}),\phi_{2}(\mu_{1},\ldots,\mu_{n})]
\]
for all $\phi_{1},\phi_{2} \in \Hom^{n}(L,L')$ and $\mu_{1},\ldots,\mu_{n} \in \Gamma L$. It also admits the following action of $G$:
\[
[R^{Hom}_{g}(\phi)](\mu_{1},\ldots,\mu_{n}) = \phi(R_{g}(\mu_{1}),\ldots,R_{g}(\mu_{n}))
\]
for all $g \in G,~\phi \in \Hom(L,L')$ and $\mu_{1},\ldots,\mu_{n} \in \Gamma L$. With this bracket and this action, $\Hom(L,L')$ becomes a PBG-LAB. At this point, let us recall the definition of an equivariant represention of a PBG-groupoid on a PBG-vector bundle (for the definition of a PBG-vector bundle see \ref{ex:frPBGgpd}). 
\begin{df}\label{df:pbggpdrep}
Let $\Upsilon \gpd P(M,G)$ be a transitive PBG-groupoid and $E \stackrel{\pi}{\rightarrow} P(M,G)$ be a PBG-vector bundle (both over the same principal bundle $P(M,G)$). An {\em equivariant representation of $\Upsilon$ on $E$} is a smooth map $\Upsilon * E \rightarrow E$ such that:
\begin{enumerate}
\item $\pi(\xi \cdot v)) = t(\xi)$ for all $(\xi , v) \in \Upsilon * E$
\item $(\xi_{1}\xi_{2}) \cdot v = \xi_{1} \cdot (\xi_{2},v)$ for all $(\xi_{1}, \xi_{2}) \in \Upsilon * \Upsilon$, $(\xi_{2}, v) \in \Upsilon * E$ 
\item $1_{\pi(v)} \cdot v = v$ for all $v \in E$.
\item $(\xi g) \cdot (v g) = (\xi \cdot v) \cdot g$.
\end{enumerate}
The set $\Upsilon * E$ consists of the pairs $(\xi,v) \in \Upsilon \times E$ such that $\alpha(\xi) = \pi(v)$.
\end{df}
Consider now the following representations of Lie groupoids:
\begin{enumerate}
\item $\Phi[L] * \Hom^{n}(L,P \times \reals) \rightarrow \Hom^{n}(L,P \times \reals)$ defined by
\[
(\xi,\phi) \mapsto \phi \circ (\xi^{-1})^{n}
\]
\item $\Phi[L] * \Hom^{n}(L,L) \rightarrow \Hom^{n}(L,L)$ defined by
\[
((\xi,\xi'),\phi) \mapsto \xi' \circ \phi \circ (\xi^{-1})^{n}
\]
\item $(\Phi[L] \times_{P \times P} \Phi[L']) * \Hom^{n}(L,L') \rightarrow \Hom^{n}(L,L')$ defined by
\[
((\xi,\xi'),\phi) \mapsto \xi' \circ \phi \circ (\xi^{-1})^{n}.
\]
\end{enumerate}
Note that the elements of the Lie groupoid $(\Phi[L] \times_{P \times P} \Phi[L]) \gpd P(M,G)$ are of the form 
\[
(\xi : E_{u} \rightarrow E_{v},\xi' : E_{u} \rightarrow E_{v})
\]
for $u,v \in P$. The Lie group $G$ acts on it by 
\[
(\xi,\xi') \cdot g = (\Hat{R}_{g} \circ \xi \circ \Hat{R}_{g^{-1}}, \Hat{R}_{g} \circ \xi' \circ \Hat{R}_{g^{-1}})
\]
for all $g \in G$, and under this action it is a PBG-groupoid. 
In \cite{LGLADG} it is proved that these representations are smooth and they induce the following representations of Lie algebroids respectively:
\begin{enumerate}
\item $\CDO[L] \rightarrow \CDO[\Hom^{n}(L,P \times \reals)]$ defined by
\[
X(\phi)(\mu_{1},\ldots,\mu_{n}) = q(X)(\phi(\mu_{1},\ldots,\mu_{n})) -
\sum_{i = 1}^{n}\phi(\mu_{1},\ldots,X(\mu_{i}),\ldots,\mu_{n})
\]
\item $\CDO[L] \rightarrow \CDO[\Hom(L,L)]$ defined by
\[
X(\phi)(\mu_{1},\ldots,\mu_{n}) = X(\phi(\mu_{1},\ldots,\mu_{n})) - \sum_{i = 1}^{n}
\phi(\mu_{1},\ldots,X(\mu_{i}),\ldots,\mu_{n})
\]
\item $\CDO[L] \oplus \CDO[L'] \rightarrow \CDO[\Hom^{n}(L,L')]$ defined by
\[
[(X \oplus X')(\phi)](\mu) = X'(\phi(\mu)) - \phi(X(\mu)).
\]
\end{enumerate}
The definition of the above representations and the fact that the standard flat connection is isometablic show that these representations are equivariant. Therefore if $\nabla, \nabla'$ are isometablic Koszul connections on $L,L'$ respectively then through the previous representations one gets isometablic Koszul connections on $\Hom^{n}(L,P \times \reals),~\Hom^{n}(L,L)$ and $\Hom(L,L')$.

\section{Equivariant transition data}

This section gives the proof of \ref{0.1} and \ref{0.2}. 

\begin{prop} \label{thm:flatisomKoszconn}
Let $L \rightarrow P(M,G)$ be a PBG-LAB. If $P(M,G)$ has a flat connection and $M$ is simply connected, then $L \Rightarrow P(M,G)$ has a flat isometablic Koszul 
connection.
\end{prop}

\pf~If $P(M,G)$ has a flat connection and $M$ is simply connected, then it is isomorphic to the trivial bundle $M \times G(M,G,pr_{1})$. So it suffices to prove 
that every PBG-LAB over a trivial principal bundle has a flat isometablic Koszul connection.

The action $\delta : (M \times G) \times G \rightarrow M \times G$ is $((x,g),h) \mapsto (x,gh)$ for all $x \in M$ and $g,h \in G$. In other words, $G$ only acts 
on itself by right translations. For every $(x,g) \in M \times G$, the partial map $\delta_{(x,g)} : G \rightarrow M \times G$ is 
$\delta_{(x,g)} = (const_{x},\ell_{g})$, where $\ell_{g} : G \rightarrow G$ denotes the left translation in $G$. Now form the action groupoid 
$(M \times G) \ract G \gpd M \times G$. This differentiates to the Lie algebroid $(M \times G) \times \gog$ with Lie bracket the one given by 
(\ref{eqn:ptwsbracket}). Here the fundamental vector field $V^{\dagger}$ corresponding to a smooth map $V : M \times G \rightarrow \gog$ is
\[
V^{\dagger}_{(x,g)} = T_{(x,e)}\delta_{(x,g)}(V_{(x,g)}) = T_{(x,e)}(const_{x},\ell_{g})
(V_{(x,g)}) = 0 \oplus T_{e}\ell_{g}(V_{(x,g)}) = T_{e}\ell_{g}(V_{(x,g)})
\]
for all $(x,g) \in M \times G$.
Since the action of $\smallskip G$ on $M \times G$ leaves $M$ unaffected, let us concentrate on the action of $G$ on itself by left translations. Form the action 
groupoid $G \ract G \gpd G.$ The map $\phi : G \ract G \rightarrow G \times G$ defined by $(g,h) \stackrel{\phi}{\mapsto} (g,gh)$
for all $g,h \in G$ makes the action groupoid isomorphic to the pair groupoid $G \times G \gpd G$. 
On the Lie algebroid level this differentiates to the isomorphism of Lie algebroids $\phi_{*} : G \ract \gog \rightarrow TG$ given by
\[
(g,V) \stackrel{\phi_{*}}{\mapsto} T_{e}\ell_{g}(V)
\]
for all $g \in G$ and $V : G \rightarrow \gog$. Since it is an isomorphism of Lie algebroids, its inverse maps the (usual) Lie bracket of $TG$ to the Lie bracket 
in $G \times \gog$ defined by (\ref{eqn:ptwsbracket}).

Now consider the derivative representation $\epsilon_{*} : (M \times G) \ract \gog \rightarrow \CDO[L]$ introduced in \ref{ex:derreps}.  Since the action of $G$ 
on $M \times G$ does not affect $M$, we can consider the Lie algebroid of the action groupoid to be $M \times (G \ract \gog)$. Finally, define 
$\nabla \co \Gamma TM \times \Gamma TG \rightarrow \Gamma \CDO[L]$ to be
\[
\nabla_{X}(V) = \epsilon_{*}(X(\phi_{*}^{-1}(V))).
\]
This is a morphism of Lie algebroids because both $\phi_{*}^{-1}$ and $\epsilon_{*}$ are. Therefore it is a flat Koszul connection. Its isometablicity is 
immediate.~\boom
\begin{prop} \label{thm:trivialisation}
If the PBG-LAB $L \rightarrow P(M,G)$ has a flat isometablic Koszul connection and $P$ is simply connected, then it is isomorphic to a trivial 
PBG-LAB (as it was introduced in example \ref{ex:trPBGLAB}).
\end{prop}

\pf~Suppose $\nabla : TP \rightarrow \CDO[L]$ is a flat isometablic Koszul connection. Then it is a morphism of PBG-algebroids, and since $P$ is simply connected 
it can be integrated to a morphism of PBG-groupoids $\phi : P \times P \rightarrow \Phi[L]$ (see \ref{prop:integrPBGmorph}). Choose a $u_{0} \in P$, define 
$\hoh = L_{u_{0}}$ and consider the following representation of $G$ on $\hoh$:
\[
\rho_{*}(g)(V) = \phi(u_{0},u_{0} \cdot g^{-1})(\Bar{R}_{g^{-1}}(V)),
\]
for all $V \in \hoh$ and $g \in G$. This map is easily shown to be a representation, to preserve the Lie bracket of $\hoh$ because 
$\Bar{R}_{g}([V,W]) = [\Bar{R}_{g}(V),\Bar{R}_{g}(V)]$ for all $V,W \in \hoh$, and $\phi(u_{0},u_{0} \cdot g^{-1}) \in \Phi[L]$. Now, the product 
$P \times \hoh$ becomes a PBG-LAB with action
\[
(u,V) \cdot g = (u \cdot g,\rho_{*}(g^{-1})(V)).
\]
Consider the map $\Psi : P \times \hoh \rightarrow L$ defined as $\Psi(u,V) = \phi(u,u_{0})(V)$. This is clearly an isomorphism of Lie algebra bundles. Moreover, 
it preserves the action because:
\begin{multline*}
\Psi((u,V) \cdot g) = \phi(u \cdot g,u_{0})(\phi(u_{0},u_{0} \cdot g)(\Bar{R}_{g}(V))) = 
\phi(u \cdot g,u_{0} \cdot g)(\Bar{R}_{g}(V)) = \\
= [\Bar{R}_{g} \circ \phi(u,u_{0}) \circ \Bar{R}_{g^{-1}} \circ \Bar{R}_{g}](V) = 
\phi(u,u_{0})(V)\cdot g = \Bar{R}_{g}(\Psi(u,V))
\end{multline*}
for all $(u,V) \in P \times \hoh$ and $g \in G$. \boom

The previous two propositions prove the existence of a special trivialization for a certain class of PBG-LABs which takes into account the group action. Namely, suppose 
$L \rightarrow P(M,G)$ is a PBG-LAB such that the Lie group $G$ is simply connected. Choose a cover $\{P_{i}\}_{i \in I}$ by principal bundle charts. That is, 
for every $i \in I$ there is an open subset $U_{i} \subseteq M$ such that $P_{i} \cong U_{i} \times G$. Without harm to the generality we may consider $U_{i}$ 
to be contractible, therefore $P_{i}$ will be simply connected. Then, \ref{thm:flatisomKoszconn} shows that the PBG-LAB $L_{P_{i}} \rightarrow P_{i}(U_{i},G)$ 
has a flat isometablic Koszul connection. Now \ref{thm:trivialisation} shows that there are Lie algebras $\hoh_{i}$ acted upon by $G$ and isomorphisms of 
PBG-LABs $\psi_{i} : P_{i} \times \hoh_{i} \rightarrow L_{P_{i}}$. If we choose one of these Lie algebras $\hoh$ and consider the composition of $\psi_{i}$ with 
a chosen isomorphism $\hoh \cong \hoh_{i}$, we obtain a trivialization for $L$ which respects the group action. 
\begin{prop}\label{LABequivsections}
Suppose $L \rightarrow P(M,G)$ is a PBG-LAB such that the Lie group $G$ is simply connected and $\{ U_{i} \}_{i \in I}$ a simple open cover of $M$. Then, for any section-atlas $\{ P_{U_{i}} \cong U_{i} \times G \}_{i \in I}$ there exists a section-atlas $\{ \Psi_{i} : P_{i} \times \hoh \rightarrow L_{P_{i}} \}_{i \in I}$ of the vector bundle $L$ such that
\[
\Psi_{i}(ug,V \cdot g) = \Psi_{i}(u,V) \cdot g
\]
for all $i \in I$, $V \in L_{U_{i}}$ and $g \in G$.
\end{prop}

Let $A \Rightarrow P(M,G)$ be a PBG-algebroid. We recall from \cite{LGLADG} the construction of the transition data $(\chi,\alpha)$. It follows from \cite[IV\S4]{LGLADG} that locally $\frac{A}{G}$ has flat connections which, due to \ref{locflatisomconn} give rise to local flat isometablic connections $\theta^{*}_{i} : TP_{i} \rightarrow A_{P_{i}}$. The transition data of $A$ is defined as:
\[
\chi_{ij} : TP_{ij} \rightarrow P_{ij} \times \hoh,~~\chi_{ij} = \Psi_{i}^{-1}(\theta^{*}_{i} - \theta^{*}_{j})
\]
and 
\[
\alpha_{ij} : P_{ij} \rightarrow Aut(\hoh),~~\alpha_{ij}(u) = \Psi_{i,u} \circ \Psi_{j,u}^{-1}
\]
Moreover, they satisfy the following, where $\Delta$ stands for the Darboux derivative:
\begin{enumerate}
\item $d\chi_{ij} + [\chi_{ij},\chi_{ij}] = 0$, i.e. each $\chi_{ij}$ is a Maurer-Cartan form,
\item $\chi_{ij} = \chi_{ik} + \alpha_{ij}(\chi_{jk})$ whenever $P_{ijk} \neq \emptyset$,
\item $\Delta(\alpha_{ij}) = \ad \circ \chi_{ij}$ for all $i, j$.
\end{enumerate}
Now \ref{0.1} follows immediately from the definition of the transition data and \ref{LABequivsections}. 

For the proof of \ref{0.2} first let us recall first the following result from \cite{Mackenzie:Liebialgds}:
\begin{thm}\label{thm:Mac-Xu}
Let $\Omega, \Xi$ be Lie groupoids over the same manifold $M$ and $\mu : A\Omega \rightarrow A\Xi$ a Lie algebroid morphism. If $\Omega$ is $s$-simply connected, then there exists a unique morphism of Lie groupoids $\phi : \Omega \rightarrow \Xi$ which differentiates to $\mu$, i.e. $\phi^{*} = \mu$.
\end{thm}
Consider a PBG-groupoid $\Upsilon \gpd P(M,G)$ and its corresponding Lie algebroid $A\Upsilon \Rightarrow P(M,G)$ with adjoint bundle $L\Upsilon$. The extension of Lie algebroids corresponding to that is
\[
\frac{L\Upsilon}{G} \inj \frac{A\Upsilon}{G} \surj \frac{TP}{G}.
\]
It follows from \cite[IV\S 4]{LGLADG} that the Lie algebroid $\frac{A\Upsilon}{G}$ (over $M$) has local flat connections $\Tilde{\theta}_{i}^{*} \co TU_{i} \rightarrow (\frac{A\Upsilon}{G})_{U_{i}}$. Due to \ref{locflatisomconn} these give rise to flat isometablic connections $\theta_{i}^{*} \co TP_{i} \rightarrow A\Upsilon_{P_{i}}$. Since the connections $\Tilde{\theta}_{i}^{*}$ are flat, they can be regarded as morphisms of Lie algebroids. With the assumption that every $U_{i}$ is contractible, and by force of the previous theorem, it follows that the $\Tilde{\theta}_{i}^{*}$s integrate uniquely to morphisms of Lie groupoids $\Tilde{\theta_{i}} : U_{i} \times U_{i} \rightarrow \frac{\Upsilon}{G}_{U_{i}}^{U_{i}}$. It was shown in the proof of \ref{basisom} that the isometablic flat connections $\theta_{i}^{*}$ corresponding to the $\Tilde{\theta}_{i}^{*}$s are in essence the maps $\tilde{\theta}_{i}^{*} \circ Tp$, therefore they also integrate uniquely to morphisms of Lie groupoids
\[
\theta_{i} \co P_{i} \times P_{i} \rightarrow \Upsilon_{P_{i}}^{P_{i}}.
\]
The uniqueness argument of \ref{thm:Mac-Xu} yields that the $\theta_{i}$s are morphisms of PBG-groupoids. That is because for every $g \in G$ the map $\theta_{i}^{g}(u,v) = \theta_{i}(ug,vg)g^{-1}$ is also a morphism of Lie groupoids and differentiates to $\theta_{i}^{*}$. It therefore follows from the uniqueness of $\theta_{i}$ that $\theta_{i}^{g} = \theta_{i}$ for all $g \in G$, consequently $\theta_{i}$ is equivariant.

Just like the non-integrable case, the aim is to show that there exists an equivariant section-atlas of the PBG-LAB $L\Upsilon$. Then, the transition data $(\chi,\alpha)$ it defines has to be equivariant. To this end, we will show that on global level there exist equivariant section-atlases for the PBG-LGB $I\Upsilon$, which differentiate to the desired atlases of $L\Upsilon$. 

First, let us give some notation. For a Lie groupoid $\Omega \rightarrow M$ we denote $\Omega_{x} = s^{-1}(\{x\})$, $\Omega^{x} = t^{-1}(\{x\})$ and $\Omega_{x}^{x}$ the Lie group $\Omega_{x} \cap \Omega^{x}$.

Now fix a basepoint $u_{0}$ in $P$ and for every $i \in I$ choose a $u_{i} \in P_{i}$ and an arrow $\xi_{i} \in \Upsilon_{u_{0}}^{u_{i}}$. Define $\sigma_{i} : P_{i} \rightarrow \Upsilon_{u_{0}}$ by
\[
\sigma_{i}(u) = \theta_{i}(u,u_{i}) \cdot \xi_{i}.
\]
This is a section of the principal bundle $\Upsilon_{u_{0}}(P,\Upsilon_{u_{0}}^{u_{0}})$. Consider the Lie group $H = \Upsilon_{u_{0}}^{u_{0}}$ and define a (left) $G$-action $\rho_{i} : G \times H \rightarrow H$ by
\[
\rho_{i}(g^{-1})(h) = \sigma_{i}(u_{i}g)^{-1} \cdot (\xi_{i}g) \cdot (hg) \cdot (\xi_{i}g)^{-1} \cdot \sigma_{i}(u_{i}g).
\]
Last, consider the sections $\psi_{i} : P_{i} \times H \rightarrow I\Omega_{P_{i}}$ of $I\Omega$ defined by
\[
\psi_{i}(u,h) = \sigma_{i}(u) \cdot h \cdot \sigma_{i}(u)^{-1}
\]
The proof of the following proposition is an immediate calculation.
\begin{prop}
The sections $\psi_{i}$ are equivariant in the sense that 
\[
\psi_{i}(ug,\rho_{i}(g^{-1})(h)) = \psi_{i}(u,h) \cdot g
\]
Their transition functions $\Tilde{\alpha}_{ij} : P_{ij} \rightarrow Aut(H)$ are equivariant in the sense
\[
\Tilde{\alpha}_{ij}(ug)(\rho_{jj}(g^{-1})(h)) = \rho_{i}(g^{-1})(\Tilde{\alpha}_{ij}(u)(h)).
\]
\end{prop}

These sections clearly differentiate to sections $\Psi : P_{i} \times \hoh \rightarrow L\Upsilon_{P_{i}}$ of $L\Upsilon$. To show that these $\Psi_{i}$s, as well as their transition functions $\Tilde{\alpha}_{ij}$ are equivariant in the sense of \ref{LABequivsections} (thus they give rise to equivariant transition data for $A\Upsilon \Rightarrow P(M,G)$), we need to show that the $G$-actions $\rho_{i}$ are local expressions of the canonical $G$-action on the PBG-LGB $I\Upsilon$. 

For every $i \in I$, consider the action groupoid $P_{i} \ract G \gpd P_{i}(U_{i},g)$ (recall example \ref{ex:actPBGgpd}) and define a map $\Tilde{\rho}_{i} \co P_{i} \ract G * I\Upsilon_{P_{i}} \rightarrow I\Upsilon_{P_{i}}$ by
\[
\Tilde{\rho}_{i}((u,g),\eta \in \Upsilon_{u}^{u}) = \psi_{i}(ug,\rho_{i}(g^{-1})(\psi_{i,u}^{-1}(\eta))).
\]
Obviously, $\pi(\Tilde{\rho}_{i}((u,g),\eta)) = ug = t(u,g)$ and $\Tilde{\rho}_{i}((u,e_{G}),\eta) = \eta$. It is a straightforward exercise to verify that
\[
\Tilde{\rho}_{i}((ug_{1},g_{2}) \cdot (u,g_{1}),\eta) = \Tilde{\rho}_{i}((ug_{1},g_{2}),\Tilde{\rho}_{i}(u,g_{1}),\eta).
\]
Also, each $\Tilde{\rho}_{i}(u,g)$ is an automorphism of $\Upsilon_{u}^{u}$, therefore it is a representation of the Lie groupoid $P_{i} \ract G$ on the Lie group bundle $I\Upsilon_{P_{i}}$, in the sense of \cite{Mac1}. The following proposition allowes us to "glue" the $\Tilde{\rho}_{i}$s together to a global map.
\begin{prop}
For all $i,j \in I$ such that $P_{ij} \neq \emptyset$, $u \in P_{ij}$, $g \in G$ and
$\eta \in \Omega_{u}^{u}$ we have
\[
\Tilde{\rho}_{i}((u,g),\eta) = \tilde{\rho}_{j}((u,g),\eta).
\]
\end{prop}

\pf The equivariance of the $\alpha_{ij}$'s gives:
\begin{multline*}
\Tilde{\rho}_{i}((u,g),\eta) =
\psi_{i}(ug,\rho_{i}(g^{-1})(\psi_{i,u}^{-1}(\eta))) =
\psi_{i}(ug,\rho_{i}(g^{-1})(\Tilde{\alpha}_{ij}(u)(\psi_{j,u}^{-1}(\eta)))) = \\
= \psi_{i}(ug,\Tilde{\alpha}_{ij}(ug)(\rho_{j}(g^{-1})(\psi_{i,u}^{-1}(\eta)))) =
\psi_{j}(ug,\rho_{j}(g^{-1})(\psi^{-1}_{i,u}(\eta))) 
= \tilde{\rho}_{j}((u,g),\eta).
\end{multline*}
\boom
Now define $\rho : (P \ract G) * I\Omega \rightarrow I\Omega$ by
$\rho((u,g),\eta \in \Omega_{u}^{u}) = \Tilde{\rho}_{i}((u,g),\eta)$, if $u \in P_{i}.$ This is a representation because each $\tilde{\rho}_{i}$ is. As a matter of fact, $\rho$ is a lot simpler than it seems. Since the sections $\{ \psi_{i} \}_{i \in I}$ are equivariant we have:
\[
\rho((u,g),\eta) = \psi_{i}(ug,\rho_{i}(g^{-1})(\psi_{i,u}^{-1}(\eta))) =
\psi_{i}(u,\psi_{i,u}^{-1}(\eta)) \cdot g = \eta \cdot g.
\]
So $\rho$ is, in fact, just the PBG structure of $I\Upsilon.$

Conversely, it is possible to retrieve the local representations $\{ \rho_{i} \}_{i \in I}$ from the PBG structure of $I\Upsilon$. Suppose $\{ \sigma_{i} : P_{i} \rightarrow \Upsilon_{u_{0}} \}_{i \in I}$ is a family of sections of $\Upsilon$. Consider the charts $\psi_{i} : P_{i}
\times H \rightarrow I\Upsilon_{P_{i}}$ defined as $\psi_{i,u}(h) =
I_{\sigma_{i}(u)}(h)$ and define $\Tilde{\rho}_{i} : P_{i} \ract G \rightarrow Aut(H)$ 
by
\[
\Tilde{\rho}_{i}(u,g)(h) = \psi_{i,ug}^{-1}(\psi_{i,u}(h) \cdot g)
\]
for all $g \in G$, $h \in H$ and $u \in P_{i}$. This is a morphism of Lie groupoids over
$P_{i} \rightarrow \cdot$. For every $i \in I$ choose $u_{i} \in P_{i}$ and define
\[
\rho_{i}(g^{-1})(h) = \Tilde{\rho}_{i}(u_{i},g)(h) = \psi^{-1}_{i,u_{i}g}
(\psi_{i,u}(h) \cdot g).
\]
Then,
\begin{multline*}
\rho_{i}(g^{-1})(h) = I_{\sigma_{i}(u_{i}g)^{-1}}(I_{\sigma_{i}(u_{i})}(h)
\cdot g) = \sigma_{i}(u_{i}g)^{-1} \cdot (\sigma_{i}(u_{i})g) \cdot (hg) \cdot
(\sigma_{i}(u_{i})^{-1}g) \cdot \sigma_{i}(u_{i}g).
\end{multline*}
The latter is exactly the original definition of the $\rho_{i}$'s. Finally, these considerations prove \ref{0.2}.

\section{Holonomy}

In this section we introduce isometablic path connections for PBG-groupoids and prove that they correspond to isometablic infinitesimal connections in the transitive case. Then we use isometablic path connections to study the holonomy of transitive PBG-groupoids. In the following we restrict to transitive PBG-groupoids $\Upsilon \gpd P(M,G)$ over a principal bundle $P(M,G)$. Let us fix some notation first. We denote $C(I,P)$ the set of continuous and piecewise smooth paths in $P$. Moreover, 
$P_{0}^{s}(\Upsilon)$ denotes the set of continuous and piecewise smooth paths $\delta : I \rightarrow \Upsilon$ which commence at an identity of $\Upsilon$ and 
$s \circ \delta : I \rightarrow M$ is constant. These sets obviously admit right actions from $G$.

\begin{df}
An {\em isometablic $C^{\infty}$ path connection} in $\Upsilon$ is a map $\Gamma : C(I,P)
\rightarrow P_{0}^{s}(\Upsilon)$, usually written $c \mapsto \bar{c}$, such that:
\begin{enumerate}
\item $\bar{c}(0) = 1_{c(0)}$ and $\bar{c}(t) \in \Upsilon^{c(t)}_{c(0)}$ for all $t \in I$.
\item If $\phi : [0,1] \rightarrow [\alpha,\beta] \subseteq [0,1]$ is a diffeomorphism then
$\overline{c \circ \phi}(t) = \bar{c}(\phi(t)) \cdot [\bar{c}(\phi(0))]^{-1}$
\item If $c \in C(I,P)$ is differentiable at $t_{0} \in I$ then $\bar{c}$ is also differentiable
at $t_{0}$.
\item If $c_{1},c_{2} \in C(I,P)$ have $\frac{dc_{1}}{dt}(t_{0}) = \frac{dc_{2}}{dt}(t_{0})$
then $\frac{d\bar{c}_{1}}{dt}(t_{0}) = \frac{d\bar{c}_{2}}{dt}(t_{0})$.
\item If $c_{1},c_{2},c_{3} \in C(I,P)$ are such that $\frac{dc_{1}}{dt}(t_{0}) +
\frac{dc_{2}}{dt}(t_{0}) = \frac{dc_{3}}{dt}(t_{0})$ then $\frac{d\bar{c}_{1}}{dt}(t_{0}) +
\frac{d\bar{c}_{2}}{dt}(t_{0}) = \frac{d\bar{c}_{3}}{dt}(t_{0})$.
\item For all $g \in G,$ $\overline{R_{g} \circ c} = \Tilde{R}_{g} \circ \bar{c}$.
\end{enumerate}
\end{df}
The first five properties in the above definition constitute the standard definition of a path
connection in a Lie groupoid, as it was postulated in \cite[II\S7]{LGLADG}. The isometablicity is
expressed by the sixth property. This definition has the following consequences which were proved in \cite[II\S7]{LGLADG}.
\begin{prop}
\begin{enumerate}
\item $\bar{\kappa}_{u} = \kappa_{1_{u}}$ for all $u \in P$ (where $\kappa_{u}$ denotes the
path constant at $u$.)
\item $\overline{c^{\leftarrow}}(t) = \bar{c}^{\leftarrow}(t) \cdot [\bar{c}(t)]^{-1}$ (where
$c^{\leftarrow}$ denotes the inverse path of $c$.)
\item $\overline{c' \cdot c} = (R_{c(1)} \circ \overline{c'}) \cdot \bar{c}$ for all composable
paths $c,c' \in C(I,P)$.
\end{enumerate}
\end{prop}
For every $c \in C(I,P)$ denote $\hat{c} \equiv c(1)$.

\begin{cor}
\begin{enumerate}
\item $\hat{\kappa}_{u} = 1_{u}$
\item $\widehat{c^{\leftarrow}} = (\hat{c})^{-1}$
\item $\widehat{c' \cdot c} = \widehat{c'} \cdot \hat{c}$.
\end{enumerate}
\end{cor}

\begin{prop}
If $c_{1},c_{2} \in C(I,P)$ are such that $\frac{dc_{1}}{dt}(t_{0}) = \lambda
\frac{dc_{2}}{dt}(t_{0})$ for some $t_{0} \in I$, then
$\frac{d\bar{c}_{1}}{dt}(t_{0}) = \lambda \cdot \frac{d\bar{c}_{2}}{dt}(t_{0})$.
\end{prop}

\begin{thm}
If $\phi : P_{U} \times (-\epsilon,\epsilon) \rightarrow P$ is a local 1-parameter group of
local transformations for $P$, where $P_{U} \cong U \times G$ is the image of a chart of the
principal bundle $P(M,G)$, then the map $\bar{\phi} \co \Upsilon^{P_{U}} \times (-\epsilon,\epsilon)
\rightarrow \Upsilon$ constructed in \cite[II\S7]{LGLADG} is a local 1-parameter group of local
transformations on $\Upsilon$ and
\begin{enumerate}
\item $t \circ \bar{\phi}_{t} = \phi_{t} \circ t$ for all $t \in (-\epsilon,\epsilon)$
\item If $\phi_{v \cdot g} = R_{g} \circ \phi_{v}$ for all $g \in G,~v \in P_{u}$ then
$\bar{\phi}_{\xi \cdot g} = \Tilde{R}_{g} \circ \bar{\phi}_{\xi}$ for all $\xi \in
\Upsilon^{P_{U}}$ and $g \in G$.
\end{enumerate}
\end{thm}

\pf~The first assertion was proved in \cite[II\S7]{LGLADG}. For the second one we have:
\[
(\Tilde{R}_{g} \circ \bar{\phi}_{\xi})(t) = \bar{\phi}_{\xi}(t) \cdot g =
\overline{R_{g} \circ \phi_{\xi}}(t) = \bar{\phi}_{\xi \cdot g}(t).
\]
\boom

\begin{thm}
There is a bijective correspondence between isometablic $C^{\infty}$ path connections
$\Gamma : c \mapsto \bar{c}$ in $\Upsilon \gpd P(M,G)$ and isometablic infinitesimal
connections $\gamma : TP \rightarrow A\Upsilon$ such that a corresponding $\Gamma$ and $\gamma$
are related by
\begin{eqnarray}\label{eqn:pathconn}
\frac{d}{dt}\bar{c}(t_{0}) = TR_{\bar{c}(t_{0})}[\gamma(\frac{d}{dt}c(t_{0}))]
\end{eqnarray}
for all $c \in C(I,P)$ and $t_{0} \in I$.
\end{thm}

\pf~Suppose given an isometablic $C^{\infty}$ path connection $\Gamma$. For $v \in P$ and
$X \in T_{v}P$ take any $c \in C(I,P)$ with $c(t_{0}) = v$ and $\frac{dc}{dt}(t_{0}) = X$ for
some $t_{0} \in I$ and define $\gamma(X) = TR_{(\bar{c}(t_{0}))^{-1}}
[\frac{d}{dt}(\bar{c}(t_{0}))]$. The smoothness and $\reals-$linearity of this connection are
proven exactly as in \cite[III\S7]{LGLADG}. We now prove its isometablicity.
\begin{multline*}
\gamma(TR_{g}(X)) = TR_{(\overline{R_{g} \circ c}(t_{0}))^{-1}}[\frac{d}{dt}
(\overline{R_{g} \circ c}(t_{0}))] = TR_{\Tilde{R}_{g} \circ (\bar{c}(t_{0}))^{-1}}
[\frac{d}{dt}(\Tilde{R}_{g} \circ \bar{c}(t_{0}))] = \\
= (TR_{(\bar{c}(t_{0}))^{-1} \cdot g} \circ T\Tilde{R}_{g})[\frac{d}{dt}(\bar{c}(t_{0}))] =  T(R_{(\bar{c}(t_{0}))^{-1} \cdot g} \circ \Tilde{R}_{g})[\frac{d}{dt}(\bar{c}(t_{0}))].
\end{multline*}
For all $\xi \in \Upsilon$ we have
\begin{multline*}
R_{(\bar{c}(t_{0}))^{-1} \cdot g} \circ \Tilde{R}_{g}(\xi) =
(\xi \cdot g) \cdot ((\bar{c}(t_{0}))^{-1} \cdot g) = (\xi \cdot (\bar{c}(t_{0}))^{-1})
\cdot g = \Tilde{R}_{g} \circ R_{(\bar{c}(t_{0}))^{-1}}(\xi).
\end{multline*}
Therefore,
\[
\ T(R_{(\bar{c}(t_{0}))^{-1} \cdot g} \circ \Tilde{R}_{g})[\frac{d}{dt}(\bar{c}(t_{0}))] =
T\Tilde{R}_{g}(TR_{(c(t_{0}))^{-1}})[\frac{d}{dt}(\bar{c}(t_{0}))] = T\Tilde{R}_{g}(\gamma(X)).
\]
Conversely, suppose given an isometablic infinitesimal connection $\gamma : TP \rightarrow A\Upsilon$.
It was shown in \cite[III\S7]{LGLADG} that given a path $c \in C(I,P)$ such that $c(0) = u_{0}$ there is a
unique path $\bar{c} \in P_{0}^{s}(\Upsilon)$ which satisfies (\ref{eqn:pathconn}). For any $g \in G$ the
path $R_{g} \circ c$ has initial data $(R_{g} \circ c)(0) = u_{0} \cdot g$. Thus there is a
unique path $\overline{R_{g} \circ c} \in P_{0}^{s}(\Upsilon)$ satisfying the differential
equation
\[
\frac{d}{dt}(\overline{R_{g} \circ c})(t_{0}) = TR_{\overline{R_{g} \circ c}(t_{0})}
[\gamma(\frac{d}{dt}(R_{g} \circ c)(t_{0}))]
\]
with initial data $\overline{R_{g} \circ c}(0) = 1_{u_{0} \cdot g}$. The path $\Tilde{R}_{g}
\circ \bar{c}$ also has initial data $\Tilde{R}_{g} \circ \bar{c}(0) = \Tilde{R}_{g}(1_{u_{0}}) =
1_{u_{0}} \cdot g = 1_{u_{0} \cdot g}$. Therefore it suffices to show that it also satisfies the
differential equation (\ref{eqn:pathconn}). Indeed:
\begin{multline*}
\frac{d}{dt}(\Tilde{R}_{g} \circ \bar{c})(t_{0}) = T\Tilde{R}_{g}(\frac{d}{dt}\bar{c}(t_{0})) =
(T\Tilde{R}_{g} \circ TR_{\bar{c}(t_{0})})[\gamma(\frac{d}{dt}c(t_{0}))] 
= T(\Tilde{R}_{g} \circ R_{\bar{c}(t_{0})})[\gamma(\frac{d}{dt}c(t_{0}))].
\end{multline*}
It is easy to see that $\Tilde{R}_{g} \circ R_{\bar{c}(t_{0})} =
R_{(\Tilde{R}_{g} \circ \bar{c})(t_{0})}$, therefore
\begin{multline*}
T(\Tilde{R}_{g} \circ R_{\bar{c}(t_{0})})[\gamma(\frac{d}{dt}c(t_{0}))] =
TR_{(\Tilde{R}_{g} \circ \bar{c})(t_{0})}[(T\Tilde{R}_{g} \circ \gamma)(\frac{d}{dt}c(t_{0}))] = \\
= TR_{(\Tilde{R}_{g} \circ c)(t_{0})}[(\gamma \circ TR_{g})(\frac{d}{dt}c(t_{0}))] =
TR_{(\Tilde{R}_{g} \circ c)(t_{0})}[\gamma(\frac{d}{dt}(R_{g} \circ c)(t_{0}))].
\end{multline*}
\boom

\begin{cor} \label{cor:exponprops}
(Of the proof.) Let $\gamma : TP \rightarrow A\Upsilon$ be an isometablic infinitesimal connection in
$\Upsilon$ and $\Gamma$ the corresponding isometablic path connection. Then:
\begin{enumerate}
\item For all $X \in \Gamma TP,$ $\Exp(t \cdot \gamma(X))(v) = \Gamma(\phi,v)(t)$ where $\phi_{t}$
is the local flow of $X$ near $v$ and $\Gamma(\phi,v) \co \reals \rightarrow \Upsilon_{v}$ is the lift
of $t \mapsto \phi_{t}(v)$.
\item The restriction of the exponential map on the image of $\gamma$ is equivariant.
\item The restriction of the exponential map on $\gamma(\Gamma^{G}TP)$ is equivariant.
\end{enumerate}
\end{cor}

\pf~The first was proved in \cite[III\S7]{LGLADG}. For the second assertion consider
$\{ \phi_{t} \}$ a local flow of $X \in \Gamma TP$ near $v \in P$. Then
$\{ \psi_{t} = R_{g} \circ \phi_{t} \circ R_{g^{-1}} \}$ is a local flow of $TR_{g}(X_{v})$.
Therefore,
\begin{multline*}
\Exp(t \cdot T\Tilde{R}_{g}(\gamma(X)))(v) = \Exp(t \cdot \gamma(TR_{g}(X)))(v)
= \Gamma(\psi,v \cdot g)(t) = \\
= \Gamma(R_{g} \circ \phi_{t} \circ R_{g^{-1}},v \cdot g)(t)
= [\Tilde{R}_{g}^{\Gamma} \circ \Gamma(\phi,v)](t) 
= \Tilde{R}_{g}^{\Gamma} \circ \Exp(t \cdot \gamma(X))(v).
\end{multline*}
The third assertion is immediate because if $X \in \Gamma^{G}TP$ then both $\{ \phi_{t} \}$ and
$\{ \psi_{t} \}$ are flows of $X$ at $v$ and $v \cdot g$ respectively. \boom

\begin{df}
Let $\Gamma$ be an isometablic $C^{\infty}$ path connection in $\Upsilon$. Then
$\Psi = \Psi(\Gamma) = \{ \bar{c}(1) : c \in C(I,P) \} \subseteq \Upsilon$ is called the
{\em holonomy subgroupoid} of $\Gamma$. The vertex $\Psi_{v}^{v}$ is the {\em holonomy group}
of $\Gamma$ at $v$.
\end{df}
We know from \cite{LGLADG} that $\Psi$ is in general an $s$-connected Lie subgroupoid of $\Upsilon
\gpd P(M,G)$. We will finish this section by showing that in case $\Gamma$ is isometablic
then it is in fact a PBG-subgroupoid of $\Upsilon$. 

\begin{prop}
Let $\Upsilon,~\Upsilon'$ be PBG-groupoids over $P(M,G)$ and $\phi \co \Upsilon \rightarrow \Upsilon'$ an
(equivariant) morphism of PBG-groupoids over $P(M,G)$. Let $\gamma$ be an isometablic
infinitesimal connection in $\Upsilon$ and $\Gamma$ its corresponding path connection. The
isometablic path connection $\Gamma'$ corresponding to the produced connection
$\gamma' = \phi_{*} \circ \gamma$ is $\Gamma' = \phi \circ \Gamma$ and its holonomy subgroupoid
is $\Psi' = \phi(\Psi)$.
\end{prop}

\pf~Take $c \in C(I,P)$. Then, since $\Bar{c}$ satisfies equation (\ref{eqn:pathconn}) for $\gamma$ it immediately follows that 
$\phi \circ \Bar{c}$ satisfies equation (\ref{eqn:pathconn}) for $\phi_{*} \circ \gamma$. Thus $\Gamma' = \phi \circ \Gamma$. \boom

\begin{df}
Let $\Upsilon,~\Upsilon'$ be PBG-groupoids over the principal bundle $P(M,G)$. The PBG-groupoid
$\Upsilon'$ is called a {\em PBG reduction} of $\Upsilon$ if there is a morphism of Lie groupoids
$\phi \co \Upsilon' \rightarrow \Upsilon$ which is equivariant and an injective immersion.
\end{df}

\begin{cor}
If $\Upsilon'$ is a PBG reduction of $\Upsilon$ and the isometablic connection
$\gamma : TP \rightarrow A \Upsilon$ takes values in $A\Upsilon'$ then $\Psi \leq \Upsilon'$
\end{cor}

\begin{thm}
Let $\Gamma$ be an isometablic path connection in $\Upsilon$. Then the holonomy subgroupoid $\Psi$
of $\Gamma$ is a PBG-subgroupoid of $\Upsilon$.
\end{thm}

\pf~The fact that $\Psi$ is a Lie subgroupoid of $\Upsilon$ is proved in \cite[III\S7]{LGLADG}.
We will only show that the action of $G$ on $\Upsilon$ can be restricted to $\Psi$. Let
$\xi = \bar{c}(1)$ for some $c \in C(I,P)$. Then $\xi \cdot g = \Tilde{R}_{g} \circ \Bar{c}(1) =
\overline{R_{g} \circ c}(1) \in \Psi$. \boom

\begin{cor} \label{cor:conninholgpd}
For each $X \in \Gamma TP$ and all $t$ sufficiently near $0$, $\gamma(X) \in \Gamma A\Psi$ and
$\Exp(t \cdot \gamma(X)) \in \Psi$.
\end{cor}

\pf~These are reformulations of \ref{cor:exponprops}. \boom

\section{Deformable sections}

The aim of this section is to prove an Ambrose-Singer theorem for transitive PBG-groupoids. The result passes through the study of deformable sections for representations of vector bundles on PBG-groupoids (see \ref{df:pbggpdrep}). Again, all groupoids regarded in this section are considered to be transitive.

\begin{df}
Let $\Upsilon \gpd M$ be a (transitive) Lie groupoid, $E \rightarrow M$ a vector bundle and
$\rho \co \Upsilon * E \rightarrow E$ a smooth representation of $\Upsilon$ on $E$. A section
$\mu \in \Gamma E$ is called {\em $\Upsilon$-deformable} if for every $x,y \in M$ there is a
$\xi \in \Upsilon^{y}_{x}$ such that $\rho(\xi,\mu(x)) = \mu(y)$.
\end{df}
If $\mu \in \Gamma E$ is $\Upsilon$-deformable then the {\em isotropy groupoid} of $\Upsilon$ at $\mu$
is $\Phi(\mu) = \{ \xi \in \Upsilon \co \xi(\mu(s(\xi))) = \mu(t(\xi)) \}$. A section $\mu$
is $\Upsilon$-deformable iff its value lies in a single orbit. The condition ensures that the
isotropy groupoid is transitive.

\begin{sloppypar}
\begin{thm} \label{thm:defsec}
Let $\Upsilon \gpd P(M,G)$ be a (transitive) PBG-groupoid and $E \rightarrow P(M,G)$ a vector bundle on
which $G$ acts by isomorphisms. Let $\rho \co \Upsilon * E \rightarrow E$ be an equivariant representation
and $\mu \in \Gamma^{G}E$. Then the following propositions are equivalent:
\begin{enumerate}
\item The section $\mu$ is $\Upsilon$-deformable.
\item The isotropy groupoid $\Phi(\mu)$ is a PBG-subgroupoid of $\Upsilon$.
\item The PBG-groupoid $\Upsilon$ has a section atlas
$\{ \sigma_{i} : P_{i} \rightarrow \Upsilon_{u_{0}} \}_{i \in I}$ such that
$\rho(\sigma_{i}(u)^{-1})(\mu(u))$ is a constant map $P_{i} \rightarrow E_{u_{0}}$.
\item The PBG-groupoid $\Upsilon$ posesses an isometablic infinitesimal connection $\gamma$ such that
\[
(\rho_{*} \circ \gamma)(\mu) = 0.
\]
\end{enumerate}
\end{thm}
\end{sloppypar}

\pf~$(({\rm i}) \Rightarrow ({\rm ii}))$. It was shown in \cite[III\S7]{LGLADG} that $\Phi(\mu)$ is a closed
and embedded Lie subgroupoid of $\Upsilon$. It remains to show that if $\xi \cdot g \in \Phi(\mu)$
for all $\xi \in \Phi(\mu)$ and $g \in G$. Indeed,
\[
\rho(\xi \cdot g,\mu(s(\xi \cdot g))) = \rho(\xi \cdot g,\mu(s(\xi)) \cdot g) =
\rho(\xi,\mu(s(\xi))) \cdot g = \mu(t(\xi)) \cdot g = \mu(t(\xi \cdot g)).
\]
\begin{sloppypar}
$(({\rm ii}) \Rightarrow ({\rm iii}))$. Since $\Phi(\mu)$ is a Lie subgroupoid of $\Upsilon$, the principal bundle
$(\Phi(\mu))_{u_{0}}(P,(\Phi(\mu))_{u_{0}}^{u_{0}},t)$ is a reduction of
$\Upsilon_{u_{0}}(P,\Upsilon_{u_{0}}^{u_{0}},t)$. Therefore there is a section atlas
$\{ \sigma_{i} : P_{U_{i}} \rightarrow \Upsilon_{u_{0}} \}_{i \in I}$ of $\Upsilon$ such that
$\sigma_{i}(u) \in (\Phi(\mu))_{u_{0}}$ for all $u \in P_{U_{i}}$. So $\sigma_{i}(u)^{-1} \in
(\Phi(\mu))_{u_{0}}$ for all $u \in P_{U_{i}}$. Equivalently, $\rho(\sigma_{i}(u)^{-1},\mu(u)) =
\mu((u_{0}))$ for all $u \in P_{i}$.
\end{sloppypar}

$(({\rm iv}) \Rightarrow ({\rm i}))$. Suppose $\Psi \leq \Upsilon$ is the holonomy PBG-subgroupoid of $\gamma$ and
$\Psi' \leq \Phi(E)$ the holonomy PBG-subgroupoid of $\rho_{*} \circ \gamma$. Then $\Psi' =
\rho(\Psi)$ and from \cite[III\S7]{LGLADG} we have:
\begin{multline*}
(\rho_{*} \circ \gamma) = 0 \Rightarrow \xi \cdot \mu(s(\xi)) = \mu(t(\xi))~
\forall \xi \in \Phi(E) \Rightarrow \rho(\eta,\mu(s(\eta))) = \mu(t(\eta))~\forall \eta \in \Upsilon.
\end{multline*}
Therefore $\mu$ is $\Upsilon$-deformable.

$(({\rm ii}) \Rightarrow ({\rm iv}))$. The isotropy subgroupoid is a PBG-groupoid, therefore it has an isometablic
connection $\gamma : TP \rightarrow A \Phi(\mu)$. This is also a connection of $\Upsilon$. From
\cite[III\S4]{LGLADG} we have that $\Gamma A\Phi(\mu) = \{ X \in \Gamma A\Upsilon \co \rho_{*}(X)(\mu) = 0\}$. 
Therefore $\rho_{*}(\gamma(X))(\mu) = 0$ for all $X \in \Gamma TP$. \boom

\begin{prop} \label{thm:defsecbasic}
Let $L$ be a vector bundle over $P(M,G)$ on which $G$ acts by isomorphisms and $[\ ,\ ]$ a section 
of the vector bundle $\Alt^{2}(L;L)$. Then the following three conditions are equivalent:
\begin{enumerate}
\item The fibers of $L$ are pairwise isomorphic as Lie algebras.
\item $L$ admits an isometablic connection $\nabla$ such that 
\[
\nabla_{X}[V,W] = [\nabla_{X}(V),W] + [V,\nabla_{X}(W)] 
\]
for all $X \in \Gamma TP$ and $V,W \in \Gamma L$.
\item $L$ is a PBG-Lie algebra bundle.
\end{enumerate}
\end{prop}

\pf~Let $\rho \co \Phi(L) * \Alt^{2}(L,L) \rightarrow \Alt^{2}(L,L)$ denote the representation defined by
\[
\rho(\xi,\phi) = \xi \circ \phi \circ (\xi^{-1} \times \xi^{-1})
\]
for all $\xi \in \Phi(L)$ and
$\phi \in \Alt^{2}(L,L)$. We have already discussed why this is an equivariant representation.
Now, $({\rm i})$ is the condition that $[\ ,\ ]$ is $\Phi[L]$-deformable and $({\rm iii})$ is the condition that $\Phi[L]$
admits a section atlas $\{ \sigma_{i} \}_{i \in I}$ such that the corresponding charts for
$\Alt^{2}(L,L)$ via $\rho$ map $[\ ,\ ] \in \Gamma Alt^{2}(L,L)$ to constant maps $P_{i}
\rightarrow \Alt^{2}(L_{u_{0}},L_{u_{0}})$. So, $({\rm i})$ and $({\rm iii})$ are equivalent by the equivalence
$(({\rm i}) \Leftrightarrow ({\rm iii}))$ of \ref{thm:defsec}. 

We also have that $\rho_{*} \co \CDO(L) \rightarrow \CDO(\Alt^{2}(L,L))$ is 
\[
\rho_{*}(D)(\phi)(V,W) = D(\phi(V,W)) - \phi(D(V),W) - \phi(V,D(W)).
\]
Therefore $({\rm ii})$ is the condition that $L$ admits an isometablic connection $\nabla$ such that
$(\rho_{*} \circ \nabla)([\ ,\ ]) = 0$. Hence $(({\rm i}) \Leftrightarrow ({\rm ii}))$ follows from the equivalence
$(({\rm i}) \Leftrightarrow ({\rm iv}))$ of \ref{thm:defsec}. \boom

\begin{prop} \label{thm:defsecvb}
Let $E$ and $E'$ be vector bundles over the principal bundle $P(M,G)$ on which $G$ acts by
isomorphisms and let $\phi : E \rightarrow E'$ be an equivariant morphism over $P$. Then the
following conditions are equivalent:
\begin{enumerate}
\item The map $P \ni u \mapsto \rk(\phi_{u}) \in \integers$ is constant.
\item There exist equivariant connections $\nabla$ and $\nabla'$ of $E$ and $E'$ respectively
such that $\nabla'_{X}(\phi(\mu)) = \phi(\nabla_{X}(\mu))$ for all $\mu \in \Gamma E$ and
$X \in \Gamma TP$.
\item There exist atlases $\{ \psi_{i} : P_{i} \times V \rightarrow E_{P_{i}} \}_{i \in I}$
and $\{ \psi'_{i} : P_{i} \times V' \rightarrow E'_{P_{i}} \}_{i \in I}$ for $E$ and $E'$
respectively such that each $\phi : E_{U_{i}} \rightarrow E'_{U_{i}}$ is of the form
$\phi_{i}(x,v) = (x,f_{i}(v))$ where $f_{i} : V \rightarrow V'$ is a linear map depending only
on $i$.
\end{enumerate}
\end{prop}

\pf~This is also an application of \ref{thm:defsec}. Consider the equivariant representation $(\Phi[E] \times_{P \times P} \Phi[E]) * \Hom(E,E') \rightarrow \Hom(E,E')$ 
we discussed in Section 3. The equivariant morphism $\phi : E \rightarrow E'$ can be regarded as an
equivariant section of the vector bundle $\Hom(E,E')$, namely assigning to every $u \in P$ the
linear map $\phi_{u} : E_{u} \rightarrow E'_{u}$. The following lemma shows that $({\rm i})$ is the
condition that $\phi$ is $(\Phi(E) \times_{P \times P} \Phi(E'))$-deformable: \\ \\
{\bf Lemma.} Let $\phi_{1} : V \rightarrow V'$ and $\phi_{2} : W \rightarrow W'$ be morphisms of
vector spaces such that $\dim V = \dim W,$ $\dim V' = \dim W'$ and $\rk(\phi_{1}) = \rk(\phi_{2})$. Then
there are isomorphisms $s : V \rightarrow W$ and $s' : V' \rightarrow W'$ such that
$s' \circ \phi_{1} = \phi_{2} \circ s$.
\begin{sloppypar}
We also discussed in section 3 that the induced representation
$\rho_{*} \co \CDO(E) \oplus \CDO(E') \rightarrow \CDO(\Hom(E,E'))$ is
$(X \oplus X')(\phi)(\mu) = X'(\phi(\mu)) - \phi(X(\mu))$. Therefore $({\rm ii})$ is exactly the
condition that $\Phi(E) \times_{P \times P} \Phi(E')$ posesses an isometablic infinitesimal
connection $\gamma$ such that $(\rho_{*} \circ \gamma)(\mu) = 0$. Finally, condition $({\rm iii})$ is
clearly condition $({\rm iii})$ of \ref{thm:defsec}. \boom
\end{sloppypar}

\begin{df}
Let $E,E'$ be two vector bundles over the principal bundle $P(M,G)$ on both of which $G$ acts by 
automorphisms. Let $\phi : E \rightarrow E'$ be an equivariant morphism of vector bundles over $M$. 
Then, $\phi$ is
\begin{enumerate}
\item {\em of locally constant rank} if $u \mapsto \rk(\phi_{u}) : P \rightarrow \integers$ is
locally constant;
\item {\em a locally constant morphism} if it satisfies condition $({\rm iii})$ of \ref{thm:defsecvb}.
\end{enumerate}
\end{df}
The proof of the following proposition is analogous to the one of \ref{thm:defsecvb}.
\begin{prop}
Let $L,L'$ be PBG-LABs on $P(M,G)$ and let $\phi : L \rightarrow L'$ be a 
morphism of PBG-LABs. Then the following conditions are equivalent:
\begin{enumerate}
\item For each $u,v \in P$ there are Lie algebra isomorphisms $\alpha : L_{u} \rightarrow L_{v}$
and $\alpha' : L'_{u} \rightarrow L'_{v}$ such that $\phi_{v} \circ \alpha = \alpha' \circ L_{u}$.
\item $L$ and $L'$ posess isometablic Lie connections $\nabla$ and $\nabla'$ respectively such
that $\phi(\nabla_{X}(V)) = \nabla'_{X}(\phi(V))$ for all $V \in \Gamma L$ and $X \in \Gamma TP$.
\item There exist Lie algebra bundle atlases $\{ \psi_{i} : P_{i} \times \gog \rightarrow
L_{P_{i}} \}_{i \in I}$ and $\{ \psi'_{i} : P_{i} \times \gog' \rightarrow
L'_{P_{i}} \}_{i \in I}$ for $L$ and $L'$ respectively such that each $\phi : L_{P_{i}}
\rightarrow L'_{P_{i}}$ is of the form $\phi(u,W) = (u,f_{i}(W))$ where $f_{i} \co \gog
\rightarrow \gog'$ is a Lie algebra morphism depending on $i$.
\end{enumerate}
\end{prop}
If $P$ is connected then $f_{i}$ may be chosen to be independent of $i$ as well.
\begin{prop} \label{thm:section}
Let $L$ be a PBG-Lie algebra bundle on $P(M,G)$ and let $L^{1}$ and $L^{2}$ be PBG-subLABs
of $L.$ Then $L^{1} \cap L^{2}$ is a PBG-subLAB of $L$ if there
is a Lie connection $\nabla$ of $L$ such that $\nabla(\Gamma L^{1}) \subseteq \Gamma L^{1}$ and
$\nabla(\Gamma L^{2}) \subseteq \Gamma L^{2}$.
\end{prop}

\pf~Without the PBG condition, this was proved in \cite[III\S7]{LGLADG}. All we need to show
is that the intersection of two PBG-subLABs is PBG, i.e. the action of $G$ on $L$
can be restricted to $L^{1} \cap L^{2}$. This is immediate. \boom

We are now ready to proceed to the proof of the Ambrose-Singer theorem for isometablic connections.
Consider a (transitive) PBG-groupoid $\Upsilon \gpd P(M,G)$ and an isometablic infinitesimal connection
$\gamma$ on $\Upsilon$.

\begin{prop} \label{prop:PBGreduct}
Let $L'$ be a PBG-subLAB of $L\Upsilon$ such that
\begin{enumerate}
\item $\Bar{R}_{\gamma}(X,Y) \in L'$ for all $X,Y \in TP$.
\item $\nabla^{\gamma}(\Gamma L') \subseteq \Gamma L'$.
\end{enumerate}
Then there is a PBG reduction $A' \leq A\Upsilon$ defined by
\[
\Gamma A' = \{ X \in \Gamma A\Upsilon : X - \gamma(q(X)) \in \Gamma L \}
\]
which has $L'$ as adjoint bundle and is such that $\gamma(X) \in A'$ for all $X \in TP$.
\end{prop}

\pf~Again, without the PBG condition this is was proved in \cite[III\S7]{LGLADG}.
We only need to show that $A'$ is a PBG-algebroid. Indeed, if $X \in \Gamma A'$ then
$\Hat{R}_{g}^{\Gamma}(X) \in \Gamma A'$ because:
\[
\Hat{R}_{g}^{\Gamma}(X) - \gamma(q(\Hat{R}_{g}^{\Gamma}(X))) =
\Hat{R}_{g}^{\Gamma}(X)(X - \gamma(q(X))) \in \Gamma L'.
\]
\boom

\begin{prop}
There is a least PBG-subLAB denoted $(L\Upsilon)^{\gamma}$ of $L\Upsilon$ which
has the properties $1$ and $2$ of \ref{prop:PBGreduct}.
\end{prop}

\pf~It suffices to prove that if $L^{1}$ and $L^{2}$ both satisfy $1$ and $2$ of \ref{prop:PBGreduct}
then $L^{1} \cap L^{2}$ does also. The only point that is not clear is that $L^{1} \cap L^{2}$
is a PBG-subLAB, and since $\nabla^{\gamma}$ is a Lie connection this is
established by \ref{thm:section}. \boom

\begin{sloppypar}
\begin{thm}
Let $\Upsilon \gpd P(M,G)$ be a PBG-groupoid and $\gamma : TP \rightarrow A\Upsilon$ an
isometablic infinitesimal connection. Let $\Gamma$ be its corresponding $C^{\infty}$ path
connection and $\Psi$ its holonomy groupoid. Then $A\Psi = (A\Upsilon)^{\gamma}$.
\end{thm}
\end{sloppypar}

\pf~We showed in \ref{cor:conninholgpd} that $\gamma$ takes values in $A\Psi$. Hence $L\Psi$ satisfies the conditions
of \ref{prop:PBGreduct} and therefore $L\Psi \geq (L\Upsilon)^{\gamma}$ and $A\Psi \geq (A\Upsilon)^{\gamma}$.
On the other hand, $\gamma$ takes values in $(A\Upsilon)^{\gamma}$ and
$A\Psi \leq (A\Upsilon)^{\gamma}$. \boom

\nocite{TFB} \nocite{Var} \nocite{Panh} \nocite{K-N:1} \nocite{K-N:2} \nocite{thesis} \nocite{Mackenzie:constronalgds} \nocite{Mackenzie:onextnsofpbs} \nocite{GHV} \nocite{Mackenzie:Liebialgds} \nocite{Dieudonne} \nocite{Moerdijk:classification}

\bibliographystyle{plain}

\end{document}